\documentclass{article}
\usepackage{graphicx} % Required for inserting images
\usepackage{amssymb}
\usepackage{latexsym}
\usepackage{url}
\usepackage{xcolor}
\definecolor{newcolor}{rgb}{.8,.349,.1}
\usepackage{algpseudocode}
\usepackage{algorithm}

\usepackage{amsthm}
\usepackage{amsfonts}
\usepackage{amsmath}
\usepackage[noblocks]{authblk}
\usepackage{graphicx}
\usepackage{geometry}
\usepackage{bm}
\usepackage{todonotes}
\usepackage{mathtools}
\usepackage{comment}
\usepackage{epstopdf}
\usepackage{float}
\usepackage{placeins}

\usepackage{xcolor}
\usepackage{graphicx}
\usepackage{geometry}
\usepackage{bm}
\usepackage{todonotes}
\usepackage{mathtools}
\newtheorem{pro}{Proposition}

\usepackage[percent]{overpic}
\usepackage[latin1]{inputenc}
\usepackage{epstopdf}
\usepackage{float}
\usepackage{placeins}
\graphicspath{ {./Figures/} }
\newcommand{\pad}[2]{\frac{\partial #1}{\partial #2}}

\title{Second order unfitted ghost-FEM for elliptic interface problems with applications to low-dimensional semiconductor devices}
\author[1]{Clarissa Astuto}
\author[2]{Giovanni Nastasi}
\affil[1]{Department of Mathematics and Computer Science, University of Catania, Catania, Italy}
\affil[2]{Department of Engineering and Architecture, University of Enna ``Kore'', Enna, Italy}
\date{}

\begin{document}

\maketitle

\begin{abstract}
We develop an unfitted ghost finite element method for elliptic interface problems with discontinuous diffusion coefficients and apply it to the electrostatic simulation of low-dimensional semiconductor devices. The proposed approach is based on a fixed Cartesian grid and represents the geometry by level-set functions, avoiding mesh generation and remeshing even in the presence of interfaces. A snapping-back-to-grid strategy is used to control the small-cut-cell issue, while interface conditions are weakly enforced by a symmetric Nitsche formulation. The method is first validated on benchmark elliptic interface problems with different geometries and coefficient jumps, showing second-order accuracy for the solution and first-order accuracy for its gradient.

As an application, we consider a graphene field-effect transistor described by a self-consistent drift-diffusion-Poisson model. The electrostatic potential is computed in a two-dimensional oxide/graphene/oxide structure, while charge transport in the graphene layer is modeled by one-dimensional bipolar drift-diffusion equations in the degenerate case and by including a field-dependent mobility model. The coupled nonlinear system is solved by a damped fixed-point iteration combined with a domain-decomposition treatment of the three-layer geometry. Numerical simulations reproduce transfer characteristics with a clear transition from an OFF state to an ON state and show the influence of the gate voltage on the two-dimensional electrostatic potential. The results also indicate that the effective thickness and discretization of the graphene layer affect the transverse potential profile, supporting the use of a full two-dimensional electrostatic description with explicit oxide/graphene interface conditions.
\end{abstract}

\section{Introduction}
{The modeling and simulation of semiconductor devices have attracted considerable attention from physicists, mathematicians, and engineers due to their fundamental role in modern technology and the challenging multiscale phenomena involved. For a comprehensive overview of the physical models and their mathematical treatment, we refer to
\cite{markowich1985stationary,markowich2012semiconductor}. }
The metal-oxide-semiconductor field-effect transistor (MOSFET) constitutes the fundamental building block of modern integrated circuits. Over several decades, extensive theoretical and numerical investigations have enabled the optimization of silicon- and III-V-semiconductor-based devices, leading to remarkable improvements in performance, scalability, and reliability \cite{selberherr,lundstrom2002fundamentals,sze2021physics}. As device dimensions continue to shrink toward the nanometer scale, however, conventional semiconductor technologies face increasing challenges associated with short-channel effects, mobility degradation, and power dissipation \cite{lundstrom2006nanoscale}. These limitations have stimulated intense research into novel materials capable of sustaining the next generation of electronic devices.

Among the emerging materials, graphene has attracted exceptional interest due to its unique physical properties. Its atomically thin structure, high carrier mobility, excellent electrical conductivity, and outstanding thermal characteristics make it a promising candidate for nanoelectronic applications \cite{CaNe,katsnelson}. In particular, the possibility of employing a conducting channel only one atomic layer thick offers an attractive route toward ultimate device scaling. Consequently, graphene field-effect transistors (GFETs) have been extensively investigated as potential alternatives or complements to conventional silicon MOSFETs.

Despite these advantages, the practical implementation of GFETs remains challenging \cite{schwierz}. The absence of an intrinsic energy bandgap in monolayer graphene results in ambipolar transport and a limited current-off region, severely restricting the achievable on/off current ratio required for digital applications. Furthermore, the transport properties of graphene supported on substrates are often degraded by impurity scattering and interactions with remote substrate phonons, leading to reduced carrier mobility. Additional complexities arise at the metal-graphene interfaces, where the contact physics differs substantially from that of conventional semiconductor devices and cannot always be described adequately through standard Ohmic or Schottky contact models.

Accurate modeling and simulation therefore play a crucial role in understanding the operation of GFETs and guiding device optimization. While drift-diffusion models remain the most widely used framework in semiconductor device simulation because of their computational efficiency, many existing GFET studies rely on simplified electrostatic descriptions, often based on reduced one-dimensional approximations of the Poisson equation \cite{ancona,jimenez,feijoo,champlain}. Such approaches may not fully capture the complex electrostatic interactions occurring in realistic device geometries. For this reason, multidimensional drift-diffusion-Poisson formulations \cite{NaRo_CNSNS} have attracted increasing attention as a means of providing a more faithful representation of charge transport and electrostatic coupling within graphene-based devices.

In \cite{NaRo_IEEE}, a novel GFET architecture was proposed that appears to mitigate the limitations associated with the zero-bandgap nature of graphene while exhibiting a sufficiently wide off-state region for practical engineering applications. The device is modeled through a self-consistent framework in which the full two-dimensional Poisson equation is coupled with the drift-diffusion transport equations. The charge in the graphene sheet is assumed to be uniformly distributed within an equivalent three-dimensional region represented by a parallelepiped whose base corresponds to the graphene surface and whose height is the sum of the thicknesses of the graphene sheet and the two graphene/oxide interfaces. In \cite{JourdanaPietra}, the electrostatic problem is reformulated as an interface Poisson problem in which the semiconducting single-layer material is modeled as a lower-dimensional manifold embedded in the computational domain, giving rise to suitable transmission conditions that account for the charge distribution and dielectric discontinuities across the interface. A comparison between the two approaches has been performed in \cite{CaNaRo_JMI}.

Another important aspect concerns the carrier statistics employed in the transport model. Owing to the high carrier densities frequently encountered in graphene, degenerate effects can become significant, making the conventional Maxwell-Boltzmann approximation potentially inadequate \cite{NaRo_JMI}. A careful assessment of the differences between degenerate and non-degenerate transport descriptions is therefore essential for obtaining reliable predictions of device performance. Moreover, the determination of accurate mobility models remains a key issue, since carrier transport in graphene is strongly influenced by multiple scattering mechanisms, including electron-phonon interactions, charged impurities, and substrate-induced effects.

In this work, we investigate the behavior of a graphene field-effect transistor through a self-consistent drift-diffusion-Poisson framework incorporating a full two-dimensional electrostatic description. Particular attention is devoted to the influence of carrier statistics and transport modeling on the predicted device characteristics. 

The discretization of the spatial derivatives is performed using the ghost-nodal finite element method (ghost-FEM), an unfitted finite element approach proposed in \cite{astuto2024nodal,astuto2024comparison}, and subsequently employed in several applications, including biological network formation \cite{astuto2023self,astuto2026modeling}, Multiscale Poisson-Nernst-Planck systems for the diffusion of ions \cite{astuto2025asymptotic,astuto2025standard}, and high-order Navier-Stokes equations \cite{dilip2026high}. Furthermore, Multigrid solvers tailored to the ghost-FEM discretization for elliptic problems on arbitrary domains have been developed and analyzed in \cite{dilip2025multigrid}. The present work focuses instead on the extension of the ghost-FEM framework to elliptic interface problems with discontinuous diffusion coefficients. Such problems introduce additional difficulties in the treatment of the solution at the interfaces where the jump conditions need to be satisfied. 

Among the numerous unfitted methods available in the literature, including finite difference schemes such as the Coco-Russo \cite{COCO2018299}, and finite element approaches such as CutFEM \cite{burman2015cutfem} or AgFEM \cite{badia2023space}, the ghost-FEM offers some advantages. In particular, it employs a very simple algorithm based on snapping-back-to-grid approach to overcome the ill-conditioning issue arising from small cut cells. Moreover, the present approach is based on a fixed Cartesian grid, consequently, neither mesh generation nor mesh adaptation to the interface is required. In standard fitted finite element frameworks, resolving complex interfaces often necessitates adaptive mesh refinement strategies (see, e.g., \cite{semplice2018adaptive}) near the interface, leading to additional meshing efforts. Even in the presence of complex or irregular interfaces, ghost-FEM employs the same squared mesh, avoiding the need for remeshing procedures. 

In the present work, ghost-FEM provides a convenient framework for coupling systems that consider different spatial dimensions, for example when considering the two-dimensional Poisson equation and the drift-diffusion equations in the one-dimensional setting in the graphene sheet. The mathematical model is defined in a semiconducting device with a central active zone between two inactive layers. To efficiently solve the resulting system, the domain is partitioned into three sub regions and a domain decomposition strategy is employed. This method provides an effective framework for the numerical solution of elliptic problems by replacing the original problem with a collection of smaller problems defined on regular subdomains \cite{DDchan1987analysis,DDchan1994domain}. A similar approach is investigated in \cite{JourdanaPietra}, for a single-layer graphene formulation. 

The paper is organized as follows. In Section \ref{SEC:math_model}, we present the drift-diffusion-Poisson mathematical model for charge transport in a GFET, together with its variational formulation. In Section \ref{sec:numerical}, we describe the numerical discretization of both the one-dimensional drift-diffusion model and the two-dimensional interface problem. In Section \ref{SEC:accuracy}, we assess the accuracy of the method in a more general framework. In Section \ref{SEC:num_res}, we present the numerical results for the GFET. Finally, Section \ref{SEC:concl} concludes the paper with final remarks and future perspectives.

\section{Mathematical model}\label{SEC:math_model}
\begin{figure}[ht]
    \centering
    \includegraphics[width=0.75\linewidth]{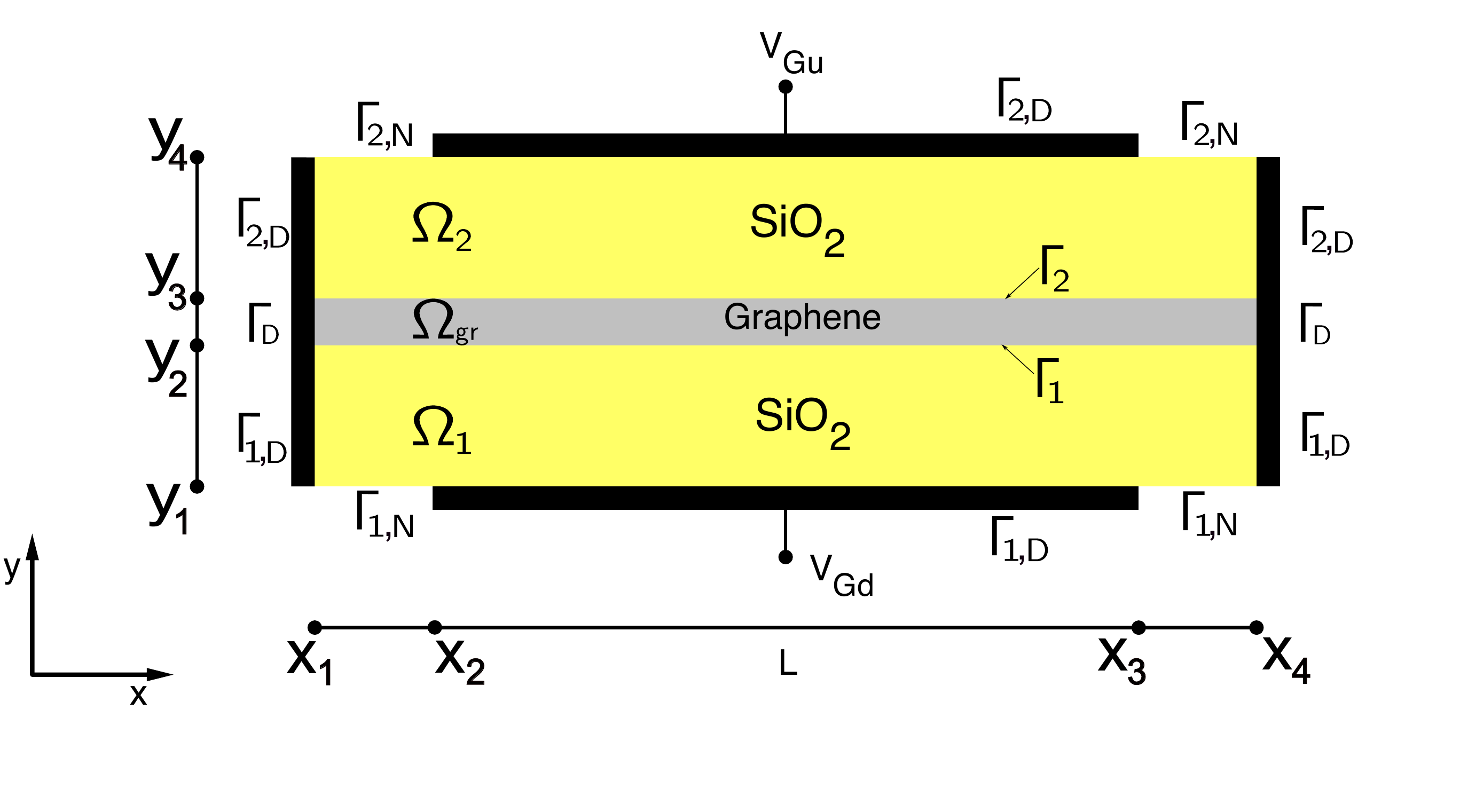}
    \caption{\textit{The domain is $\Omega = \Omega_1\cup\Omega_{\rm gr}\cup\Omega_2$ and its boundary $\partial \Omega$ is composed by the union of Dirichlet and Neumann boundaries, $\Gamma_D,\Gamma_{ 1,D}, \Gamma_{ 2,D}$ and $\Gamma_{1, N}, \Gamma_{2, N}$ respectively. The two interfaces are $\Gamma_{1}$ and $\Gamma_{2}$.}}
    \label{fig:domain}
\end{figure}

A GFET typically consists of a single-layer (or few-layer) graphene sheet acting as the channel, connected to source and drain electrodes. In a back-gate GFET, the channel is electrostatically controlled by a gate electrode located beneath the substrate and separated from graphene by a dielectric layer. In a double-gate GFET, an additional top gate is placed above the graphene channel, providing enhanced electrostatic control. The combined action of the top and bottom gates enables more efficient carrier modulation and improved device performance. Fig.~\ref{fig:domain} shows the cross-sectional view of the GFET structure introduced in \cite{NaRo_IEEE}. The mathematical model is defined in the rectangular domain $\Omega=[x_1,x_4]\times[y_1,y_4]$, which is partitioned into three subdomains. The regions $\Omega_1=[x_1,x_4]\times[y_1,y_2]$ and $\Omega_2=[x_1,x_4]\times[y_3,y_4]$ represent the lower and upper SiO$_2$ dielectric layers, respectively, while $\Omega_{\rm gr}=[x_1,x_4]\times[y_2,y_3]$ denotes the graphene layer. Although graphene is only one atom thick and charge transport occurs exclusively along the $x$-direction, a finite thickness $t_{\rm gr}=y_3-y_2$ is introduced to accurately describe the electrostatic potential. This thickness accounts for the charge distribution across the graphene layer and the adjacent oxide-graphene interfaces, following the approach adopted in \cite{NaRo_CNSNS,NaRo_IEEE}. The metallic source, drain, and gate contacts (represented in black) are not explicitly included in the computational domain; their effect is incorporated through suitable boundary conditions imposed on the corresponding portions of the domain boundary.

The charge transport is described by the one-dimensional bipolar stationary drift-diffusion equations, which read
\begin{subequations}
\label{eq:currents}
    \begin{alignat}{2}
        \pad{J_n}{x}  & = e(G-R) && \qquad {\rm in } \, \Omega_{\rm gr}\\
        \pad{J_p}{x}  & = -e(G-R) && \qquad {\rm in } \, \Omega_{\rm gr}
    \end{alignat}
\end{subequations}
where $J_n$ and $J_p$ are the electron and hole current density, respectively, and $G-R$ denotes the generation-recombination term. The positive elementary charge is denoted by $e$. In the non-degenerate case, the current densities are expressed in terms of the quasi-Fermi energies $\varepsilon^{(n)}$ and $\varepsilon^{(p)}$ by
\begin{subequations}
\label{eq:currentsnp}
    \begin{align}
        J_n & = \mu_n n\pad{\varepsilon^{(n)}}{x} \\
         J_p & = \mu_p p \pad{\varepsilon^{(p)}}{x} 
    \end{align}
\end{subequations}
where $n$ and $p$ represent the electron and hole densities, respectively, and $\mu_n$ and $\mu_p$ are appropriate mobility models. The carrier densities are related to the quasi-Fermi energy by
\begin{subequations}
\label{eq:np_relations}
    \begin{align}
        n(\varepsilon^{(n)}-\varepsilon_D) & =\int_0^{+\infty} f_{\mathrm{FD}}(\varepsilon;\varepsilon^{(n)}-\varepsilon_D) \, d\varepsilon, \\
         p(\varepsilon^{(p)}-\varepsilon_D) & = \int_0^{+\infty} f_{\mathrm{FD}}(\varepsilon;-(\varepsilon^{(p)}-\varepsilon_D)) \, d\varepsilon,
    \end{align}
\end{subequations}
where $f_{\mathrm{FD}}$ is the Fermi-Dirac statistics
\begin{equation}
    \label{FD}
    f_{\mathrm{FD}}(\varepsilon;\varepsilon_F-\varepsilon_D) = \frac{D(\varepsilon)}{1+\exp\left(\frac{\varepsilon-(\varepsilon_F-\varepsilon_D)}{k_B T} \right)}.
\end{equation}
In \eqref{FD}, $k_B$ is the Boltzmann constant, $T$ is the temperature (kept constant at 300 K), $\varepsilon_D$ is the Dirac-point energy, and $D(\varepsilon)$ is the density of states in graphene, i.e.,
\begin{equation}
    D(\varepsilon) = \frac{2}{\pi}\frac{\varepsilon}{(\hbar v_F)^2},
\end{equation}
with $\hbar$ the reduced Planck constant and $v_F$ the Fermi velocity. In graphene, the physically meaningful reference energy is the Dirac point, where the conduction and valence bands touch. The electrostatic potential $\varphi$ shifts the whole band structure in energy, so the Dirac point itself moves. Therefore, we have $\varepsilon_D=-e\varphi$.

When graphene is connected to metallic source and drain contacts, the metal and graphene must reach electrochemical equilibrium. Because the metal and graphene generally have different work functions, electrons are transferred across the metal/graphene interface until their Fermi levels align. For copper contacts, this charge transfer shifts the graphene Fermi level relative to the Dirac point. In the quoted case, Cu induces an $n$-type doping, meaning that the Fermi level lies above the Dirac point, i.e.
\begin{equation}
    \varepsilon_F-\varepsilon_D = \Delta\varepsilon >0.
\end{equation}
In our simulations, we will adopt $\Delta\varepsilon = 0.25$ eV \cite{NaRo_IEEE}. Therefore, we impose Dirichlet boundary conditions
\begin{subequations}
    \begin{align}
        \varepsilon_l &:= \varepsilon^{(n)}(x = x_1) = \varepsilon^{(p)}(x = x_1) = \Delta\varepsilon -e\varphi(x_1,y), \qquad \forall y\in[y_2,y_3],\\
        \varepsilon_r &:= \varepsilon^{(n)}(x = x_4) = \varepsilon^{(p)}(x = x_4) = \Delta\varepsilon -e\varphi(x_4,y), \qquad \forall y\in[y_2,y_3].
    \end{align}
\end{subequations}

We remark that, although $J_n$ and $J_p$ depend on both $x$ and $y$, charge transport in graphene occurs only along the $x$-direction. {Consequently, in Section~\ref{sec:VF}, we define the one-dimensional function space $V = H^1([x_1,x_4])$ to introduce the variational formulation of the drift-diffusion equations.} Therefore, the total current density flowing through the contact, including both electron and hole contributions, is obtained by averaging the current densities across the effective graphene thickness $t_{\rm gr}$ (see Fig.~\ref{fig:charact}):
\begin{equation}
\label{eq:Jav}
J^{\mathrm{av}}(x) = \frac{1}{t_{\rm gr}} \int_{y_2}^{y_3} \left( J_n(x,y)+J_p(x,y) \right) \,dy .
\end{equation}

In our simulations, we adopt the mobility model proposed in \cite{Dorgan}, that is
\begin{equation}
    \mu_n(n,E) = \frac{\nu_s(n)}{\left( 1+\left( \frac{\nu_s E}{v_{\mathrm{sat}}} \right)^\zeta \right)^{1/\zeta}},
\end{equation}
where $E=\left\vert -\frac{\partial\varphi}{\partial x} \right\vert$, $v_{\mathrm{sat}}$ is the saturation velocity (we take the value 0.2 $\mu$m/ps), $\zeta=2$ and
\begin{equation}
    \nu_s(n) = \frac{\mu_0}{\left( 1+\frac{n}{n_{\mathrm{ref}}} \right)^\sigma},
\end{equation}
where $\mu_0$ is the low field mobility. The values of  $\mu_0,n_\mathrm{ref} $ and $\sigma$ are defined in Table~\ref{table_parameters}. The hole mobility $\mu_p$ takes the same expression as $\mu_n$. 

In the considered operating regime, the device geometry suppresses the activation of minority carriers, so that transport is mainly governed by the self-consistent electrostatic potential. Moreover, since a reliable spatially resolved generation--recombination model for graphene is still debated and other previous GFET simulations \cite{NaRo_CNSNS,NaRo_IEEE,NaRo_JMI} showed only a weak sensitivity to the recombination time, we neglect this contribution and set \(G-R=0\).

The drift-diffusion equations are coupled with the Poisson equation for the electrostatic potential $\varphi \in \Omega$,
\begin{equation}
\label{eq:PoissonPhi}
    -\nabla \cdot (\epsilon \nabla \varphi) =\gamma \qquad {\rm in} \, \Omega
\end{equation}
where 
\begin{align}
\label{eq:l_expr}
   \gamma = \left\{\begin{array}{cc}
       e(p - n + n_{\rm imp})/t_{\rm gr}  & \qquad {\rm in} \, \Omega_{\rm gr} \\
       0  & \qquad {\rm otherwise} 
    \end{array} \right.
\end{align}
and $n_{\rm imp}$ is the areal density of the impurity charges at the graphene/oxide interface. The dielectric coefficient $\epsilon$ is given by
\begin{align}
    \epsilon = \left\{\begin{array}{cc}
       \epsilon_{\rm gr} & \qquad {\rm in} \, \Omega_{\rm gr}  \\
       \epsilon_{\rm ox}  & \qquad {\rm otherwise}
    \end{array} \right.
\end{align}
with $\epsilon_{\rm gr}$ and $\epsilon_{\rm ox}$ dielectric constants in graphene and oxide, respectively.

Eq.~\eqref{eq:PoissonPhi} is augmented with the following boundary conditions in $\Gamma = \partial \Omega$ for $\varphi$
\begin{subequations} \label{eq:BCphi}
    \begin{align}
        \varphi &= 0 \quad {\rm at} \, x = x_1,\\
        \varphi &= V_{\rm b} \quad {\rm at} \, x = x_4,\\ 
        \varphi & = V_{\rm G_d}  \quad {\rm at} \, y = y_1, \, x\in [x_2,x_3],  \\ 
        \varphi & = V_{\rm G_u}  \quad {\rm at} \, y = y_4, \, x\in [x_2,x_3],  \\
        \nabla_{\bf n} \varphi & = 0 \quad {\rm at} \, y = y_1, y_4, \, x < x_2 \, \lor x > x_3, 
    \end{align}
\end{subequations}
where $V_{\rm b}$ is the bias voltage, $V_{\rm G_d}$ is the down gate-source potential and  $V_{\rm G_u}$ is the upper gate-source potential. $V_{\rm G_d}$ and $V_{\rm G_u}$ include the (subtracted) flat-band voltages. All the values of the parameters are shown in Table~\ref{table_parameters}.

Finally, we define the interface conditions
\begin{subequations} \label{eq:BCphi2}
    \begin{align}
    \label{eq:bcsolu}
    [|\varphi|] & = 0 \quad {\rm in} \, \Gamma_{\rm 1}\cup\Gamma_{\rm 2}  \\ \label{eq:bcder} [|\epsilon \nabla_{\bf n} \varphi|] & = 0 \quad {\rm in} \, \Gamma_{\rm 1}\cup\Gamma_{\rm 2} 
    \end{align}
\end{subequations}
where $\nabla_{\bf n}$ denotes the outgoing normal derivative. For the sake of simplicity, we refer to Fig.~\ref{fig:domain} for the description of the setting.
The domain $\Omega$ is partitioned into three subdomains, so that
$\Omega=\Omega_1\cup\Omega_{\rm gr}\cup\Omega_2$.
Moreover, we decompose the boundary of $\Omega_1$ as
$\partial\Omega_1=\Gamma_1\cup\Gamma_{1,D}\cup\Gamma_{1,N}$,
where $\Gamma_{1,D}$ and $\Gamma_{1,N}$ denote the portions of the boundary where Dirichlet and Neumann boundary conditions are imposed, respectively. The boundary $\Gamma_1$ represents the interface between $\Omega_1$ and $\Omega_{\rm gr}$. Analogously, the boundary of $\Omega_2$ is decomposed as
$\partial\Omega_2=\Gamma_2\cup\Gamma_{2,D}\cup\Gamma_{2,N}$, where $\Gamma_2$ denotes the interface between $\Omega_2$ and $\Omega_{\rm gr}$.

\subsection{Variational formulation}
\label{sec:VF}
In this section, we derive the variational formulation of the coupled problem defined by systems~(\ref{eq:currents}-\ref{eq:currentsnp}) and (\ref{eq:PoissonPhi}-\ref{eq:l_expr}), together with the associated boundary \eqref{eq:BCphi} and interface conditions \eqref{eq:BCphi2}.

To this end, we consider the one-dimensional drift-diffusion model in (\ref{eq:currents}-\ref{eq:currentsnp}) and define $V=H^1([x_1,x_4])$.
Multiplying Eqs.~(\ref{eq:currents}-\ref{eq:currentsnp}) by a test function $\psi\in V$ and integrating over the interval $[x_1,x_4]$, we obtain
\begin{subequations}
\label{eq:variational}
    \begin{align}
       \left.\left(\mu_n  n\pad{\varepsilon^{(n)}}{x}\psi(x)\right)\right|_{x_1}^{x_4} + {\left.\left(\mu_n  n\varepsilon^{(n)} \pad{\psi}{x}\right)\right|_{x_1}^{x_4}} - \int_{x_1}^{x_4} \mu_n n\pad{\varepsilon^{(n)}}{x} \pad{\psi}{x} \, dx + \lambda \left( \varepsilon^{(n)}(x_1)\psi(x_1) + \varepsilon^{(n)}(x_4)\psi(x_4) \right) + & \\= \int_{x_1}^{x_4}  f\psi(x) \, dx + {\left(\mu_n(x_1) n(x_1)\varepsilon_l \psi'(x_1) + \mu_n(x_4) n(x_4)\varepsilon_r \psi'(x_4)\right)} + \lambda \left( \varepsilon_l\psi(x_1) + \varepsilon_r\psi(x_4) \right) , \quad \forall \psi\in V \\ 
        \left.\left(\mu_p  p\pad{\varepsilon^{(p)}}{x}\psi(x)\right)\right|_{x_1}^{x_4} + {\left.\left(\mu_p p\varepsilon^{(p)} \pad{\psi}{x}\right)\right|_{x_1}^{x_4}}-  \int_{x_1}^{x_4} \mu_p  p\pad{\varepsilon^{(p)}}{x} \pad{\psi}{x} \, dx + \lambda \left( \varepsilon^{(p)}(x_1)\psi(x_1) + \varepsilon^{(p)}(x_4)\psi(x_4) \right) + & \\= - \int_{x_1}^{x_4}  f\psi(x) \, dx +{\left(\mu_p(x_1) p(x_1)\varepsilon_l \psi'(x_1) + \mu_p(x_4) p(x_4)\varepsilon_r \psi'(x_4)\right)} + \lambda \left( \varepsilon_l\psi(x_1) + \varepsilon_r\psi(x_4) \right) , \quad \forall \psi\in V
    \end{align}
\end{subequations}
with $f = e(G-R).$ We apply Nitsche technique to symmetrize the variational formulation and add a penalization term $\lambda\to\infty$, that will force Dirichlet conditions; see
\cite{Courant,HilbertCourant}. The choice of $\lambda$ is discussed in Section~\ref{sec:snap}.

Now we look at the system~(\ref{eq:PoissonPhi}-\ref{eq:l_expr}) with boundary conditions in \eqref{eq:BCphi} and interface conditions in \eqref{eq:BCphi2}. We define the space of test functions $W = H^1(\Omega)$, multiply Eq.~\eqref{eq:PoissonPhi} by $\psi \in W$, and integrate over the domain $\Omega$, obtaining
\begin{align}
\label{eq:phi_var}
    -\int_\Omega \nabla \cdot \left( \epsilon \nabla \varphi \right) \psi \,d\Omega = \int_\Omega \gamma\,\psi  \,d\Omega.
\end{align}
Considering the contribution of each subdomain, Eq.~\eqref{eq:phi_var} becomes 
\begin{align}
\label{eq:1variational2D}
    -\epsilon_{\rm ox}\int_{\Omega_1} \Delta \varphi \, \psi \,d\Omega -\epsilon_{\rm ox}\int_{\Omega_2} \Delta \varphi \, \psi \,d\Omega -\epsilon_{\rm gr}\int_{\Omega_{\rm gr}} \Delta \varphi \, \psi \,d\Omega = \int_\Omega \gamma \psi \,d\Omega.
\end{align}
Integrating by parts, we obtain the following
\begin{subequations}
\label{eq:2variational2D}
\begin{align}
    \epsilon_{\rm ox}\int_{\Omega_1} \nabla \varphi \cdot \nabla \psi \,d\Omega - \epsilon_{\rm ox}\int_{\partial \Omega_1} \pad{\varphi}{n} \psi \,dl + \epsilon_{\rm ox}\int_{\Omega_2} \nabla \varphi \cdot \nabla \psi \,d\Omega - \epsilon_{\rm ox}\int_{\partial \Omega_2} \pad{\varphi}{n} \psi \,dl + \\
+ \epsilon_{\rm gr}\int_{\Omega_{\rm gr}} \nabla \varphi \cdot \nabla \psi \,d\Omega - \epsilon_{\rm gr}\int_{\partial \Omega_{\rm gr}} \pad{\varphi}{n} \psi \,dl  =  \int_\Omega \gamma \psi \,d\Omega.
\end{align}
\end{subequations}
Here we make use of the boundary conditions in \eqref{eq:BCphi} and of the interface conditions in \eqref{eq:BCphi2}, and Eq.~\eqref{eq:2variational2D} becomes
\begin{subequations}
\label{eq:3variational2D_all}
\begin{align}
\label{eq:3variational2D}
\epsilon_{\rm ox}\int_{\Omega_1} \nabla \varphi \cdot \nabla \psi \,d\Omega - \epsilon_{\rm ox}\int_{\Gamma_{1,D}} \pad{\varphi}{n} \psi \,dl - \epsilon_{\rm ox}\int_{\Gamma_1} \pad{\varphi}{n} \psi \,dl + \\
    + \epsilon_{\rm ox}\int_{\Omega_2} \nabla \varphi \cdot \nabla \psi \,d\Omega - \epsilon_{\rm ox}\int_{\Gamma_{2,D}} \pad{\varphi}{n} \psi \,dl - \epsilon_{\rm ox}\int_{\Gamma_2} \pad{\varphi}{n} \psi \,dl + \\
+ \epsilon_{\rm gr}\int_{\Omega_{\rm gr}} \nabla \varphi \cdot \nabla \psi \,d\Omega - \epsilon_{\rm gr}\int_{\Gamma_D} \pad{\varphi}{n} \psi \,dl - \epsilon_{\rm gr}\int_{\Gamma_1} \pad{\varphi}{n} \psi \,dl - \epsilon_{\rm gr}\int_{\Gamma_2} \pad{\varphi}{n} \psi \,dl =  \int_\Omega \gamma \psi \,d\Omega.
\end{align}
\end{subequations}
To make symmetric the variational formulation, we apply Nitsche technique first to $\Gamma_{1,D}$ and $\Gamma_{2,D}$. %and leave the treatment of the two interfaces $\Gamma_1$ and $\Gamma_2$ to Section~\ref{appendix_steps} to incorporate the jump conditions. 
Thus, Eq.~\eqref{eq:3variational2D_all} becomes
\begin{subequations}
%\label{eq:5variational2D}
\begin{align*}
\epsilon_{\rm ox}\int_{\Omega_1} \nabla \varphi \cdot \nabla \psi \,d\Omega - \epsilon_{\rm ox}\int_{\Gamma_{1,D}} \left( \pad{\varphi}{n} \psi + \varphi \pad{\psi}{n} \right) \,dl - \epsilon_{\rm ox}\int_{\Gamma_1} \pad{\varphi}{n} \psi \,dl + \\
    + \epsilon_{\rm ox}\int_{\Omega_2} \nabla \varphi \cdot \nabla \psi \,d\Omega - \epsilon_{\rm ox}\int_{\Gamma_{2,D}} \left( \pad{\varphi}{n} \psi + \varphi \pad{\psi}{n} \right) \,dl - \epsilon_{\rm ox}\int_{\Gamma_2} \pad{\varphi}{n} \psi \,dl + \\
+ \epsilon_{\rm gr}\int_{\Omega_{\rm gr}} \nabla \varphi \cdot \nabla \psi \,d\Omega - \epsilon_{\rm gr}\int_{\Gamma_D} \left( \pad{\varphi}{n} \psi + \varphi \pad{\psi}{n} \right) \,dl - \epsilon_{\rm gr}\int_{\Gamma_1} \pad{\varphi}{n} \psi \,dl - \epsilon_{\rm gr}\int_{\Gamma_2} \pad{\varphi}{n} \psi \,dl = \\
\int_\Omega \gamma \psi \,d\Omega - \epsilon_{\rm ox}\int_{\Gamma_{1,D}}  \left.\varphi \right|_{\Gamma_{1,D}}\pad{\psi}{n} \,dl - \epsilon_{\rm ox}\int_{\Gamma_{2,D}}  \left.\varphi \right|_{\Gamma_{2,D}}\pad{\psi}{n} \,dl.   
\end{align*}
\end{subequations}
Now, we add the penalization terms to impose the Dirichlet boundary conditions
\begin{subequations}
%\label{eq:5variational2D}
\begin{align*}
\epsilon_{\rm ox}\int_{\Omega_1} \nabla \varphi \cdot \nabla \psi \,d\Omega - \epsilon_{\rm ox}\int_{\Gamma_{1,D}} \left( \pad{\varphi}{n} \psi + \varphi \pad{\psi}{n} \right) \,dl + \lambda \int_{\Gamma_{1,D}} \varphi\, \psi  \,dl - \epsilon_{\rm ox}\int_{\Gamma_1} \pad{\varphi}{n} \psi \,dl + \\
    + \epsilon_{\rm ox}\int_{\Omega_2} \nabla \varphi \cdot \nabla \psi \,d\Omega - \epsilon_{\rm ox}\int_{\Gamma_{2,D}} \left( \pad{\varphi}{n} \psi + \varphi \pad{\psi}{n} \right) \,dl + \lambda \int_{\Gamma_{2,D}} \varphi\, \psi  \,dl - \epsilon_{\rm ox}\int_{\Gamma_2} \pad{\varphi}{n} \psi \,dl + \\
+ \epsilon_{\rm gr}\int_{\Omega_{\rm gr}} \nabla \varphi \cdot \nabla \psi \,d\Omega - \epsilon_{\rm gr}\int_{\Gamma_D} \left( \pad{\varphi}{n} \psi + \varphi \pad{\psi}{n} \right) \,dl  + \lambda \int_{\Gamma_{D}} \varphi\, \psi  \,dl  - \epsilon_{\rm gr}\int_{\Gamma_1} \pad{\varphi}{n} \psi \,dl - \epsilon_{\rm gr}\int_{\Gamma_2} \pad{\varphi}{n} \psi \,dl = \\
\int_\Omega \gamma \psi \,d\Omega - \epsilon_{\rm ox}\int_{\Gamma_{1,D}}  \left.\varphi \right|_{\Gamma_{1,D}}\pad{\psi}{n} \,dl - \epsilon_{\rm ox}\int_{\Gamma_{2,D}}  \left.\varphi \right|_{\Gamma_{2,D}}\pad{\psi}{n} \,dl + \lambda \int_{\Gamma_{1,D}} \left.\varphi \right|_{\Gamma_{1,D}}\, \psi  \,dl + \lambda \int_{\Gamma_{2,D}} \left.\varphi \right|_{\Gamma_{2,D}} + \lambda \int_{\Gamma_{D}} \left.\varphi \right|_{\Gamma_{D}}\, \psi \, dl.    
\end{align*}
\end{subequations}
Finally, we make use of the interface conditions \eqref{eq:BCphi2} in $\Gamma_1$ and $\Gamma_2$
\begin{subequations}
\label{eq:4variational2D}
\begin{align}
\label{eq:der1}
\epsilon_{\rm ox}\int_{\Omega_1} \nabla \varphi \cdot \nabla \psi \,d\Omega - \epsilon_{\rm ox}\int_{\Gamma_{1,D}} \left( \pad{\varphi}{n} \psi + \varphi \pad{\psi}{n} \right) \,dl + \lambda \int_{\Gamma_{1,D}} \varphi\, \psi  \,dl - \epsilon_{\rm gr}\int_{\Gamma_1} \pad{\varphi}{n} \psi \,dl + \\ \label{eq:jder2}
    + \epsilon_{\rm ox}\int_{\Omega_2} \nabla \varphi \cdot \nabla \psi \,d\Omega - \epsilon_{\rm ox}\int_{\Gamma_{2,D}} \left( \pad{\varphi}{n} \psi + \varphi \pad{\psi}{n} \right) \,dl + \lambda \int_{\Gamma_{2,D}} \varphi\, \psi  \,dl - \epsilon_{\rm gr}\int_{\Gamma_2} \pad{\varphi}{n} \psi \,dl + \\ \label{eq:solu}
+ \epsilon_{\rm gr}\int_{\Omega_{\rm gr}} \nabla \varphi \cdot \nabla \psi \,d\Omega - \epsilon_{\rm gr}\int_{\Gamma_D} \left( \pad{\varphi}{n} \psi + \varphi \pad{\psi}{n} \right) \,dl  + \lambda \int_{\Gamma_{D}} \varphi\, \psi  \,dl  - \epsilon_{\rm gr}\int_{\Gamma_1 \cup \Gamma_2} \left( \pad{\varphi}{n} \psi + \varphi \pad{\psi}{n} \right)  \,dl = \\ % - \epsilon_{\rm gr}\int_{\Gamma_2} \pad{\varphi}{n} \psi \,dl = \\
\int_\Omega \gamma \psi \,d\Omega - \epsilon_{\rm ox}\int_{\Gamma_{1,D}}  \left.\varphi \right|_{\Gamma_{1,D}}\pad{\psi}{n} \,dl - \epsilon_{\rm ox}\int_{\Gamma_{2,D}}  \left.\varphi \right|_{\Gamma_{2,D}}\pad{\psi}{n} \,dl + \\ +\lambda \int_{\Gamma_{1,D}} \left.\varphi \right|_{\Gamma_{1,D}}\, \psi  \,dl + \lambda \int_{\Gamma_{2,D}} \left.\varphi \right|_{\Gamma_{2,D}} \, \psi \, dl  + \lambda \int_{\Gamma_{D}} \left.\varphi \right|_{\Gamma_{D}}\, \psi \, dl - \epsilon_{\rm gr}\int_{\Gamma_1 \cup \Gamma_2} \left.\varphi \right|_{\Gamma_{1}\cup\Gamma_2}\pad{\psi}{n} dl,  
\label{eq:phiDir}
\end{align}
\end{subequations}
where, in Eqs.~(\ref{eq:der1}-\ref{eq:jder2}), we make use of the interface conditions in \eqref{eq:bcder}, and in Eq.~\eqref{eq:solu} of the condition \eqref{eq:bcsolu}. In other words, in the latter one, we symmetrize the term $\epsilon_{\rm gr}\int_{\Gamma_1 \cup \Gamma_2} \left( \pad{\varphi}{n} \psi + \varphi \pad{\psi}{n} \right)  \,dl$ in \eqref{eq:solu}, and evaluate $ \epsilon_{\rm gr}\int_{\Gamma_1 \cup \Gamma_2} \left.\varphi \right|_{\Gamma_{1}\cup\Gamma_2}\pad{\psi}{n} dl$ on the right hand side of the equation, in \eqref{eq:phiDir}. We leave the description of the choice of the parameter $\lambda$ to Section~\ref{sec:snap}.
% ,  to incorporate the jump conditions. 
% The parameter $\lambda$ is chosen in such a way that $\lambda = h^{-\alpha}$, where $\alpha$ is defined in Section~\ref{sec:snap}.

The numerical scheme used to solve the two variational formulations is described in Section~\ref{sec:numerical}. The detailed treatment of the interfaces $\Gamma_1$ and $\Gamma_2$ is deferred to Section~\ref{appendix_steps}. Since the coupling between the three subdomains is entirely mediated through these interfaces, their formulation plays a central role in the numerical method. A dedicated section allows us to describe the numerical computation of the coupling conditions and the resulting decomposition strategy.

\section{Numerical scheme}
\label{sec:numerical}
In this section, we adopt a one- and two-dimensional space discretizations based on finite elements method for arbitrary domains. The two-dimensional domain is $\Omega = [x_1,x_4]\times[y_1,y_4] \subset S$, with  $S = [0,L_S]^2$ a square region. The side length $L_S$ of $S$ is chosen such that $L_S>L$, ensuring sufficient space for the definition and labeling of ghost points when Dirichlet or Neumann boundary conditions are imposed. In our simulations $L_S = 0.102.$ Moreover, the domain is partitioned as $\Omega = \Omega_1\cup\Omega_{\rm gr} \cup \Omega_2$, where $\Omega_1 = [x_1,x_4]\times[y_1,y_2], \, \Omega_{\rm gr} = [x_1,x_4]\times[y_2,y_3], $ and $\Omega_2 = [x_1,x_4]\times[y_3,y_4]$ (see Fig.~\ref{fig:domain}). The subdomain $\Omega_{\rm gr}$ is treated as a collection of one-dimensional domains (horizontal lines), on which the numerical scheme is applied to solve the corresponding one-dimensional problem (see Fig.~\ref{fig:levelset}). Finally, we set by $L = x_4-x_1$ and $H = y_4 - y_1$, the two space lengths. 

\subsection{Space discretization of the one-dimensional drift-diffusion model}
For the system~(\ref{eq:currents}-\ref{eq:currentsnp}) in $\Omega_{\rm gr}$, we consider the space $V_h$ of dimension $N+1$ as the space of the splines $S_1$
on the mesh of size $h = L/N$, $N\in \mathbb N$, where $N$ is the number of subintervals that define the mesh 
\[
\mathcal T_h=\{K_i=[\xi_i,\xi_{i+1}]\}_{i=1}^{N}, \quad \xi_1 = x_1, \quad \xi_{N+1} = x_4.
\]
As basis for the space $V_h$ we choose continuous, piecewise polynomial functions defined on the interval $[x_1,x_4]$. 

Since $\Omega_{\rm gr} = [x_1,x_4]\times[y_2,y_3]$, there exists an integer $N_{\rm gr} \in \mathbb{N}$ denoting the number of grid points in the $y$-direction within $\Omega_{\rm gr}$. We compute the solution of the one-dimensional problem \eqref{eq:variational} for each $y^j_{\rm gr},\, j = 1,\cdots,N_{\rm gr}$. 
To emphasize this dependence, we denote by $f_{j,h} \approx f(x_h,y^j_{\rm gr}), \, \, j = 1,\cdots,N_{\rm gr}$, and apply the same argument for all the variables in Eq.~\eqref{eq:variational_Vh}. Specifically, in Eq.~\eqref{eq:variational}, the functions $\varepsilon^{(n)}, \varepsilon^{(p)}, \mu_n, \mu_p, n, p, f \in V$ are approximated by $\varepsilon_h^{(n)}, \varepsilon_h^{(p)}, (\mu_n n)_{j,h}, (\mu_p p)_{j,h}, f_{j,h} \in V_h$, respectively.  
To perform computations, the domain $\Omega_{\rm gr}$ is approximated by a set of $N_{\rm gr}$ one-dimensional domains (see Fig.~\ref{fig:levelset}). Consequently, the original integrals, defined over $\Omega_{\rm gr}$, are now evaluated as follows.
%%%%%%%%%%%%%%%%%%%%%%%%%%%%%%%%%%%%%%%%%%%%%%%%%%%%%%%%%%%5
\begin{pro}\label{pro:variational_2}
Given $f_{j,h} \in L^2(\mathcal T_h)$ for all $j = 1, \cdots, N_{\rm gr}$ find   $\varepsilon^{(n)}_h, \varepsilon^{(p)}_h \in V_h$ such that it holds
\begin{subequations}
\label{eq:variational_Vh}
    \begin{align}
       \left.\left((\mu_n n)_{j,h}\pad{\varepsilon_h^{(n)}}{x}\psi_h\right)\right|_{x_1}^{x_4} + {\left.\left((\mu_n  n)_{j,h}\varepsilon_h^{(n)} \pad{\psi_h}{x}\right)\right|_{x_1}^{x_4}} -\left( (\mu_n n)_{j,h}\pad{\varepsilon_h^{(n)}}{x}, \pad{\psi_h}{x} \right) + \lambda \left( \varepsilon_h^{(n)}(x_1)\psi_h(x_1) + \varepsilon_h^{(n)}(x_4)\psi_h(x_4) \right) + & \\= \left( f_{j,h},\psi_h \right) + {\left((\mu_n n)_{j,h}(x_1)\varepsilon_l \psi'(x_1) + (\mu_n n)_{j,h}(x_4)\varepsilon_r \psi'(x_4)\right)} + \lambda \left( \varepsilon_l\psi_h(x_1) + \varepsilon_r\psi_h(x_4) \right) , \qquad \forall \psi_h \in V_h,\\ 
        \left.\left((\mu_p p)_{j,h}\pad{\varepsilon_h^{(p)}}{x}\psi_h\right)\right|_{x_1}^{x_4} + {\left.\left((\mu_p  p)_{j,h}\varepsilon_h^{(p)} \pad{\psi_h}{x}\right)\right|_{x_1}^{x_4}} -\left( (\mu_p p)_{j,h}\pad{\varepsilon_h^{(p)}}{x}, \pad{\psi_h}{x} \right) + \lambda \left( \varepsilon_h^{(p)}(x_1)\psi_h(x_1) + \varepsilon_h^{(p)}(x_4)\psi_h(x_4) \right) + & \\= -\left( f_{j,h},\psi_h \right) + {\left((\mu_p p)_{j,h}(x_1)\varepsilon_l \psi'(x_1) + (\mu_p p)_{j,h}(x_4)\varepsilon_r \psi'(x_4)\right)} + \lambda \left( \varepsilon_l\psi_h(x_1) + \varepsilon_r\psi_h(x_4) \right), \qquad \forall \psi_h \in V_h.
        \end{align}
\end{subequations}
and $ \forall j = 1, \cdots, N_{\rm gr}$, where $(\cdot,\cdot)$ denotes the scalar product in ${L^2(\mathcal T_h)}$. The parameter $\lambda$ is chosen in such a way that $\lambda = h^{-\alpha}$, where $\alpha$ is defined in Section~\ref{sec:snap}.
\end{pro}
\subsection{Space discretization of the interface problem in two dimensions}
The square region $S$ is discretized by a uniform mesh of squares ({\it cells}) of size $h^2, \,h=L/N$, where $N \in \mathbb N$ is the number of cells in $x$ and $y$ direction, respectively.  
The set of cells is denoted by $\mathcal C$, with $\# \mathcal C = N^2$. A cell $K \in \mathcal C$ is classified as {\it active cell} if $K \cap  \Omega \neq \emptyset$, or {\it inactive cell} otherwise.
%The union of all active cells is denoted by $\Omega_\text{act}$.
An active cell $K$ can be {\it internal cell} if $K \subseteq \Omega$ or {\it cut cell} otherwise (namely $K \cap \Gamma \neq \emptyset)$. 

Following the approach shown in  \cite{Osher2002,Russo2000,book:72748, Sussman1994}, the domains $\Omega_i,\, i = 1,2,{\rm gr},$ are implicitly defined by a level set function $\phi_i(x,y)$ that is negative inside $\Omega_i$, positive in $S\setminus \Omega_i$ and zero on the boundary $\partial \Omega_i$:
\begin{equation}
	\Omega_i = \{(x,y): \phi_i(x,y) < 0\}, \qquad
	\partial \Omega_i = \{(x,y): \phi_i(x,y) = 0\}, \quad i = 1,2,{\rm gr}.
\end{equation}

The set of grid points is denoted by $\mathcal{N}$, with cardinality $\# \mathcal{N} = (N+1)^2$. It is partitioned as $\mathcal N = \mathcal I \cup \mathcal G \cup \mathcal O$, corresponding to the  internal, ghost and inactive points, respectively. The internal points are defined as those grid points for which the level-set function takes negative values.
While, the set of ghost points, $\mathcal{G}$, consists of grid points that belong to the complement domain $\Omega_c \coloneqq S \setminus \Omega$ { and are vertices of cut cells}, formally defined as
\begin{equation}
\notag
	(x,y) \in \mathcal{G} \iff (x,y) \in {\mathcal N}\cap \Omega_c  \text{ and } \{(x \pm h,y),(x,y\pm h), (x \pm h,y\pm h) \} \cap \mathcal I \neq \emptyset.
\end{equation}
All the other points, are inactive points and belong to $\mathcal O$.

Lastly, we define the computational domain $\Omega^h = \Omega \cap \mathcal C$  (or $\Omega_1^h = \Omega_1\cap \mathcal C$, $\Omega_2^h = \Omega_2\cap \mathcal C$ and $\Omega_{\rm gr}^h = \Omega_{\rm gr}\cap \mathcal C$). In each cut cell, the portion of boundary is approximated by a straight line. The union of these segments is denoted by $\Gamma^h = \Gamma^h_{1} \cup \Gamma^h_2 \cup \Gamma^h_{1,D} \cup \Gamma^h_{2,D} \cup \Gamma^h_{1,N} \cup \Gamma^h_{2,N}$ (see Fig.~\ref{fig:domain}).
Consequently, the integrals in \eqref{eq:4variational2D} are now evaluated over $\Omega^h$ and $\Gamma^{h}$, respectively.

The discrete space ${W}_h$ is constructed from continuous, piecewise polynomial functions defined on $\Omega^h$.
The functions $\psi$, which are defined in the function space $ W$, are approximated by the function $\psi_h$ that belong to finite-dimensional subspace ${W}_h$, defined as follows:
\begin{align}\label{eq:Vh}
{W}_h & = \{ {\psi}_h \in {W} : {\psi}_h|_K \in \mathbb Q_1(K), \quad \forall K\in \Omega^h\}, 
\end{align}
where $\mathbb Q_1(K)$ denotes the space of bilinear Lagrangian shape functions on the element $K$.
%%%%%
% The discrete space $W_h$ is given by the piecewise bilinear functions which are continuous in $S$.
% As a basis of $W_h$, we choose the following functions:
% \begin{equation}
%     \psi_{i}(x,y) = \max\left\{
%         \left(1-\frac{|x-x_i|}{h}\right),0
%     \right\}\max\left\{\left(1-\frac{|y-y_i|}{h}\right),0
%     \right\},
%     \label{eq:V_h2}
% \end{equation}
% with $i = (i_1,i_2)$ an index that identifies a node on the grid. The generic element $w_h\in W_h$ has the following representation
% \begin{equation}
%     w_h(x,y) = \sum_{i\in\nodes}w_i \psi_i(x,y).
%     \label{eq:u_h2}
% \end{equation}

As in the one-dimensional case, to obtain a discrete variational problem, we consider \eqref{eq:4variational2D} with both trial and test functions taken in $W_h$, leading to the following discrete problem. 
\begin{pro}\label{pro:variational_2D}
Given $\gamma_h \in L^2(\Omega^h)$, find 
 $\varphi_h \in W_h$ such that it holds 
\begin{subequations}
\label{eq:variational2D_Wh}
    \begin{align}
\epsilon_{\rm ox}\left( \nabla \varphi_h, \nabla \psi_h \right)_{L^2(\Omega^h_1)} - \epsilon_{\rm ox}\left( \pad{\varphi_h}{n}, \psi_h \right)_{L^2(\Gamma^h_{1,D})} - \epsilon_{\rm ox} \left(\varphi_h, \pad{\psi_h}{n} \right)_{L^2(\Gamma^h_{1,D})} + \lambda\left( \varphi_h, \psi_h \right)_{L^2(\Gamma^h_{1,D})}- \epsilon_{\rm gr}\left(\pad{\varphi_h}{n}, \psi_h \right)_{L^2(\Gamma^h_1)} + \\
    + \epsilon_{\rm ox}\left( \nabla \varphi_h, \nabla \psi_h \right)_{L^2(\Omega^h_2)} - \epsilon_{\rm ox} \left( \pad{\varphi_h}{n}, \psi_h\right)_{L^2(\Gamma^h_{2,D})} - \epsilon_{\rm ox} \left(\varphi_h, \pad{\psi_h}{n} \right)_{L^2(\Gamma^h_{2,D})} + \lambda\left( \varphi_h, \psi_h \right)_{L^2(\Gamma^h_{2,D})}- \epsilon_{\rm gr}\left( \pad{\varphi_h}{n}, \psi_h\right)_{L^2(\Gamma^h_2)} + \\
+ \epsilon_{\rm gr}\left( \nabla \varphi_h, \nabla \psi_h\right)_{L^2(\Omega_{\rm gr})} - \epsilon_{\rm gr}\left( \pad{\varphi_h}{n}, \psi_h \right)_{L^2(\Gamma^h_D)}  -\epsilon_{\rm gr} \left(\varphi_h, \pad{\psi_h}{n} \right)_{L^2(\Gamma^h_D)} - \epsilon_{\rm gr}\left( \pad{\varphi_h}{n}, \psi_h\right)_{L^2(\Gamma^h_1)}  - \epsilon_{\rm gr}\left( {\varphi_h}, \pad{\psi_h}{n}\right)_{L^2(\Gamma^h_1)} + \\ +\lambda\left( \varphi_h, \psi_h \right)_{L^2(\Gamma^h_{D})}  - \epsilon_{\rm gr}\left( \pad{\varphi_h}{n}, \psi_h\right)_{L^2(\Gamma^h_2)} - \epsilon_{\rm gr}\left( {\varphi_h}, \pad{\psi_h}{n}\right)_{L^2(\Gamma^h_2)} 
= \\
\left( \gamma_h, \psi_h \right)_{L^2(\Omega^h)} - \epsilon_{\rm ox}\left(  \left.\varphi_h \right|_{\Gamma_{1,D}},\pad{\psi_h}{n} \right)_{L^2(\Gamma^h_{1,D})} - \epsilon_{\rm ox}\left(  \left.\varphi_h \right|_{\Gamma_{2,D}},\pad{\psi_h}{n} +  \right)_{L^2(\Gamma^h_{2,D})} +  \lambda\left( \left.\varphi_h \right|_{\Gamma_{1,D}}, \psi_h \right)_{L^2(\Gamma^h_{1,D})} + \\+ \lambda\left( \left.\varphi_h \right|_{\Gamma_{2,D}}, \psi_h \right)_{L^2(\Gamma^h_{2,D})} + \lambda\left( \left.\varphi_h \right|_{\Gamma_{D}}, \psi_h \right)_{L^2(\Gamma^h_{D})} - \epsilon_{\rm gr}\left( \left.{\varphi_h}\right|_{\Gamma^h_1}, \pad{\psi_h}{n}\right)_{L^2(\Gamma^h_1)} - \epsilon_{\rm gr}\left( \left.{\varphi_h}\right|_{\Gamma^h_2}, \pad{\psi_h}{n}\right)_{L^2(\Gamma^h_2)} , \quad \forall \psi_h \in W_h. 
    \end{align}
\end{subequations}
\end{pro}

\subsection{Snapping-back-to-grid technique for interface problems}
\label{sec:snap}
\begin{figure}
    \centering
\begin{minipage}[b]
		{.49\textwidth}
		\centering
	\begin{overpic}[abs,width=\textwidth,unit=1mm,scale=.25]{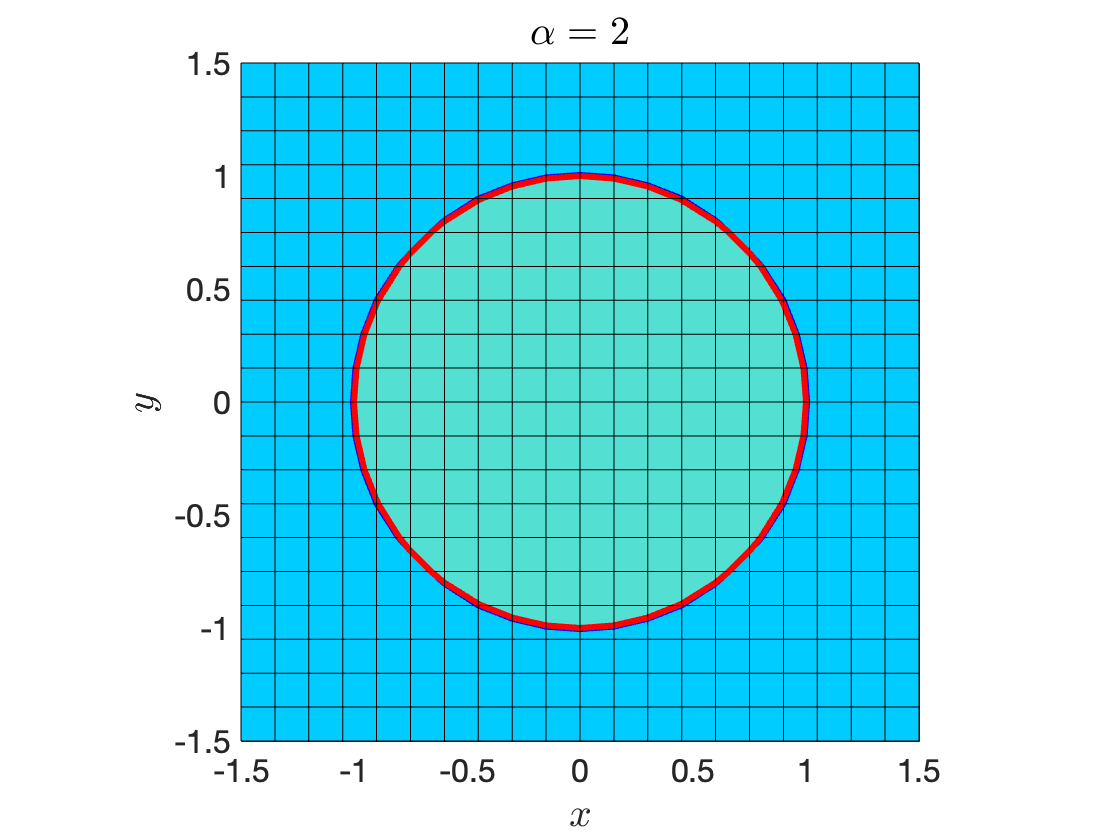}
    \put(30,55){\LARGE $\Omega_1$}
    \put(35,30){\LARGE $\Omega_2$}
\end{overpic}
\end{minipage}         
\begin{minipage}[b]
		{.49\textwidth}
		\centering
	\begin{overpic}[abs,width=\textwidth,unit=1mm,scale=.25]{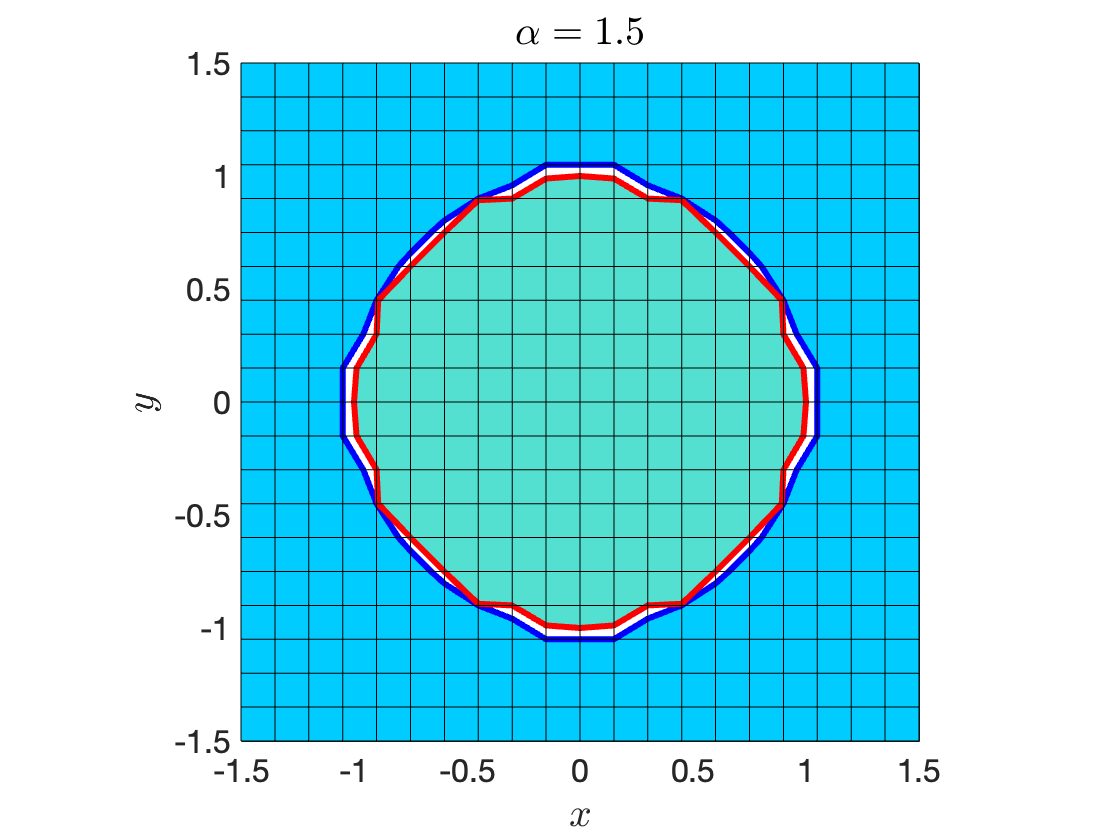}
\put(30,55){\LARGE $\Omega_1$}
\put(35,30){\LARGE $\Omega_2$}
\end{overpic}
\end{minipage}          
    \caption{\textit{Representation of the domains after the snapping-back-to-grid procedure (Algorithm \ref{alg_snap}), with $\alpha = 2$ in (a) and $\alpha = 1.5$ in (b).}}
    \label{fig:snap}
\end{figure}

\begin{figure}
    \centering
\begin{minipage}[b]
		{.49\textwidth}
		\centering
	\begin{overpic}[abs,width=\textwidth,unit=1mm,scale=.25]{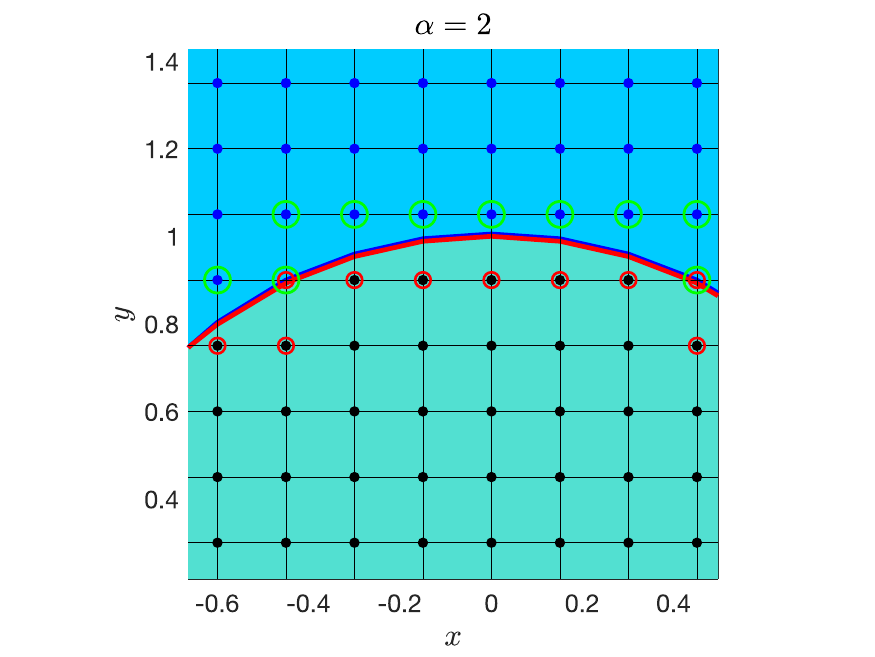}
\put(30,55){\LARGE $\Omega^h_1$}
\put(30,26){\LARGE $\Omega^h_2$}
\end{overpic}
\end{minipage}         
\begin{minipage}[b]
		{.49\textwidth}
		\centering
	\begin{overpic}[abs,width=\textwidth,unit=1mm,scale=.25]{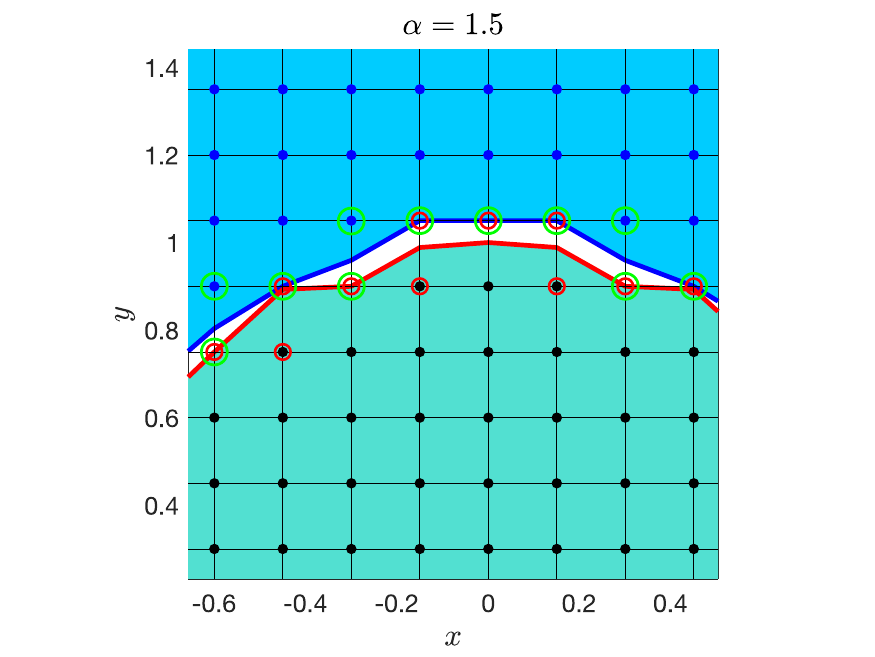}
\put(30,54.5){\LARGE $\Omega^h_1$}
\put(30,25.5){\LARGE $\Omega^h_2$}
\end{overpic}
\end{minipage}          
    \caption{\textit{Points classification after the snapping-back-to-grid procedure (Algorithm \ref{alg_snap}) for $\alpha = 2$ (a) and $\alpha = 1.5$ (b). Blue and black dots represent interior points of $\Omega^h_1$ and $\Omega^h_2$, respectively. Red and green circles denote ghost points associated with $\Omega^h_1$ and $\Omega^h_2$, respectively.}}
    \label{fig:snapls}
\end{figure}

When dealing with arbitrary domains, the so-called \textit{small cut} issue arises, which might result in losing the coercivity of the bilinear form. We avoid this problem by perturbing the mesh geometry.
In particular, in \cite{astuto2024nodal}, the perturbation procedure we introduce is well known to the engineering community by the name of snapping-back-to-grid. In other words, whenever the distance between an interior point and the boundary is smaller than $h^\alpha$, the corresponding cell is discarded. Fig.~\ref{fig:snap}-\ref{fig:snapls} illustrate the effect of this procedure for $\alpha=2$ (left panel) and $\alpha=1.5$ (right panel). As expected, smaller values of $\alpha$ lead to the removal of a larger number of cells. Consequently, the computational domain becomes closer to a polygonal approximation of the original geometry, and a larger number of points change their classification from interior to ghost points (see Fig.~\ref{fig:snapls}).
This technique improves the conditioning of the resulting linear systems by eliminating arbitrarily small cut cells. On the other hand, the geometric approximation error increases and the second-order accuracy is no longer guaranteed, as illustrated by the results reported in Table~\ref{tab:accuracy}. 
The parameter $\alpha$ is chosen based on the accuracy. Since the numerical scheme is second-order accurate, choosing $\alpha>2$ we do not expect any further improvement in the convergence rate. Therefore, we restrict our investigation to $\alpha=1.5$, $1.75$, and $2$.

\begin{algorithm}[H]
\caption{Snapping back to grid}\label{alg_snap}
\begin{algorithmic}
\For{$k \in \mathcal N$} 
\If{$\phi(k)<0 \quad \& \quad |\phi(k)|<h^{ \alpha}$} 
   \State $\phi(k) := {\tt eps} $ 
\EndIf 
\EndFor
\end{algorithmic}
\end{algorithm}

\section{Accuracy tests}
\label{SEC:accuracy}
In this section, we present different numerical tests of the one-dimensional drift-diffusion problem and of the two dimensional elliptic interface problem. Our main goal is to examine the accuracy and convergence behavior of the numerical scheme. We test different geometries, and for each one we show the accuracy of the solution, which is computed as follows:
\begin{subequations}
\label{eq:error}
\begin{align}
    {\rm error}_i & = \frac{\sqrt{\int_{\Omega_i} |u_h - u^i_{\rm exa}|^2 dx dy}}{\sqrt{\int_{\Omega_i} |u^i_{\rm exa}|^2 dx dy}}, \quad i = 1,2\\ 
    \nabla{\rm error}_i & = \frac{\sqrt{\int_{\Omega_i} |\nabla u_h - \nabla u^i_{\rm exa}|^2 dx dy}}{\sqrt{\int_{\Omega_i} |\nabla u^i_{\rm exa}|^2 dx dy}}, \quad i = 1,2,
\end{align}
\end{subequations}
where the integral over $\Omega_i$ is computed using the midpoint quadrature rule.
\subsection{Drift-diffusion model in one dimension}
Here, we perform numerical tests to check that the numerical discretization for the drift-diffusion equation achieves the expected accuracy. We consider the following problem, whose exact solution is known. 
\begin{subequations}
\label{eq:DD_1D}
\begin{align}
    \pad{}{x}\left( \nu(x)\pad{u}{x}\right) & = f(x) \quad [0,1]\\
    u(0) & = 0, \qquad u(1) = 0
\end{align}
\end{subequations}
with $\nu(x) = \exp(\cos(x))$ and $f(x) = -2\pi\exp(\cos(x))(\sin(x)\cos(2\pi x) + 2\pi\sin(2\pi x))$, whose exact solution has the following expression
\begin{equation}
	u_{\rm exa} = \sin(2\pi x).
	\label{eq:uexa1D}
\end{equation}

\begin{figure}
    \centering
\begin{minipage}[b]
		{.49\textwidth}
		\centering
	\begin{overpic}[abs,width=\textwidth,unit=1mm,scale=.25]{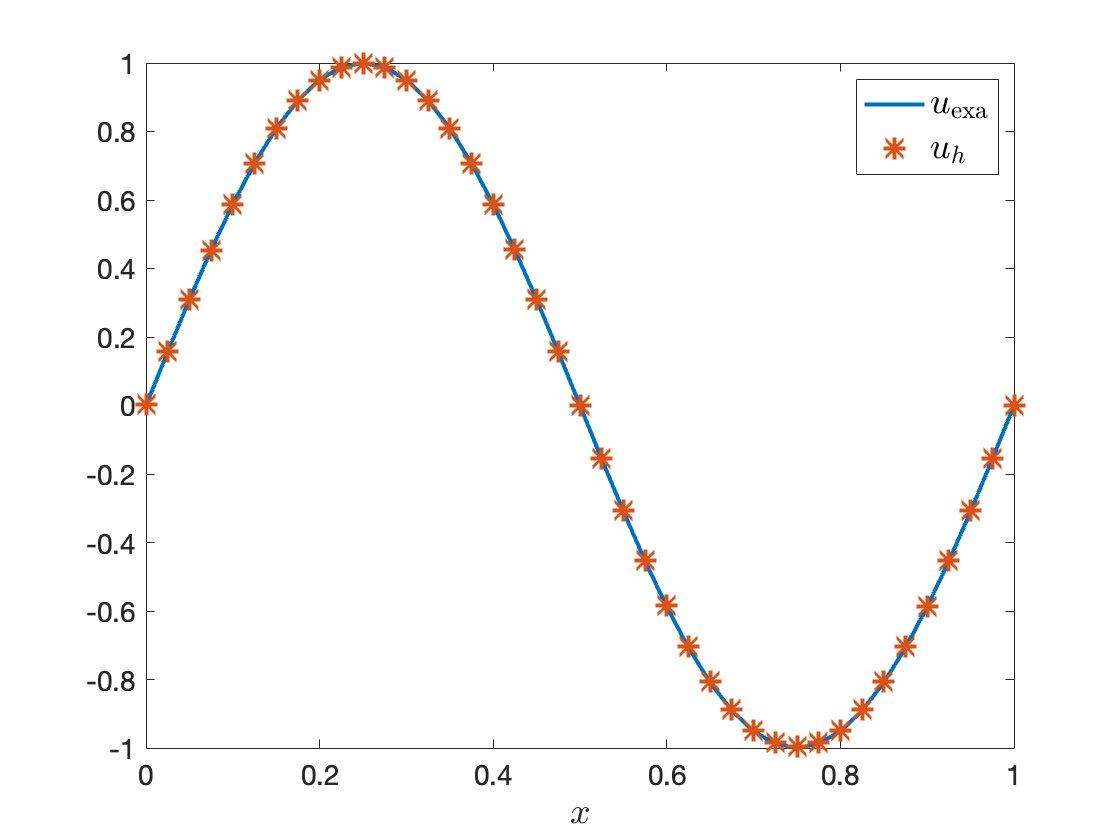}
    \put(2,65){(a)}
\end{overpic}
\end{minipage}         
\begin{minipage}[b]
		{.49\textwidth}
		\centering
	\begin{overpic}[abs,width=\textwidth,unit=1mm,scale=.25]{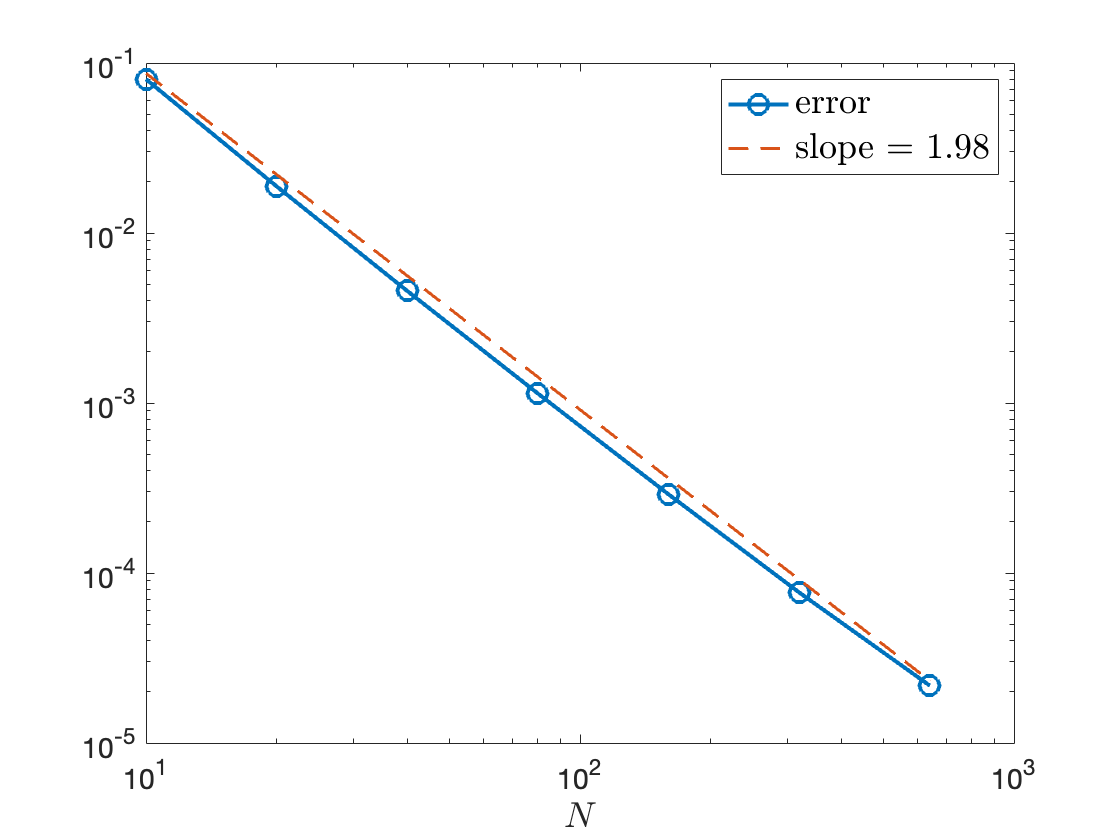}
    \put(2,65){(b)}
\end{overpic}
\end{minipage}          
    \caption{\textit{Panel (a) compares the exact solution (blue line), given by Eq.~\eqref{eq:uexa1D}, with the numerical solution (red asterisks) for the problem defined in Eq.~\eqref{eq:DD_1D}, with number of elements $N = 40$. Panel (b) shows the convergence order of the numerical scheme.}}
\label{fig:accuracy_1D}
\end{figure}

In Fig.~\ref{fig:accuracy_1D} (panel (a)), we show the comparison between the exact solution (blue line) and the numerical solution (red asterisks) with number of elements $N = 40$. In panel (b), we show the convergence order as a function of $N$.

\subsection{Interface problem in two dimensions}
\label{sec:2D_ns}
\begin{figure}
    \centering   
\begin{minipage}[b]
		{.32\textwidth}
		\centering
	\begin{overpic}[abs,width=\textwidth,unit=1mm,scale=.25]{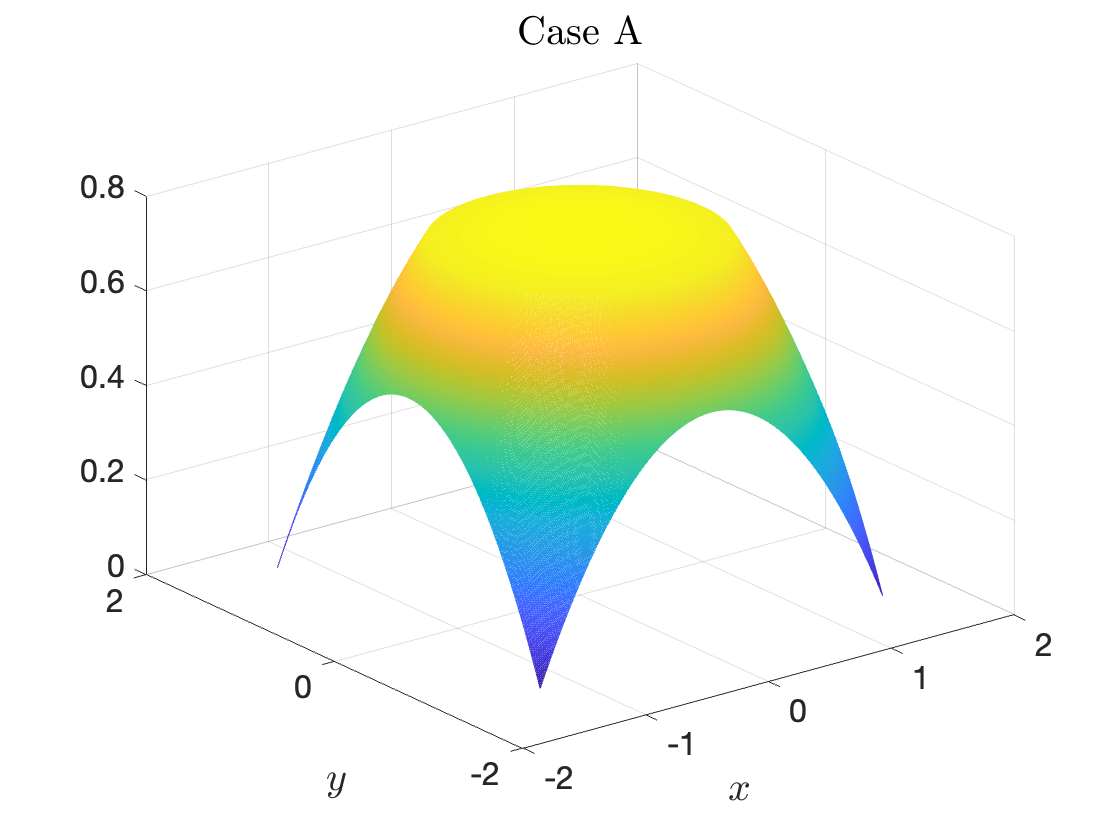}
\end{overpic}
\end{minipage}         
\begin{minipage}[b]
		{.32\textwidth}
		\centering
	\begin{overpic}[abs,width=\textwidth,unit=1mm,scale=.25]{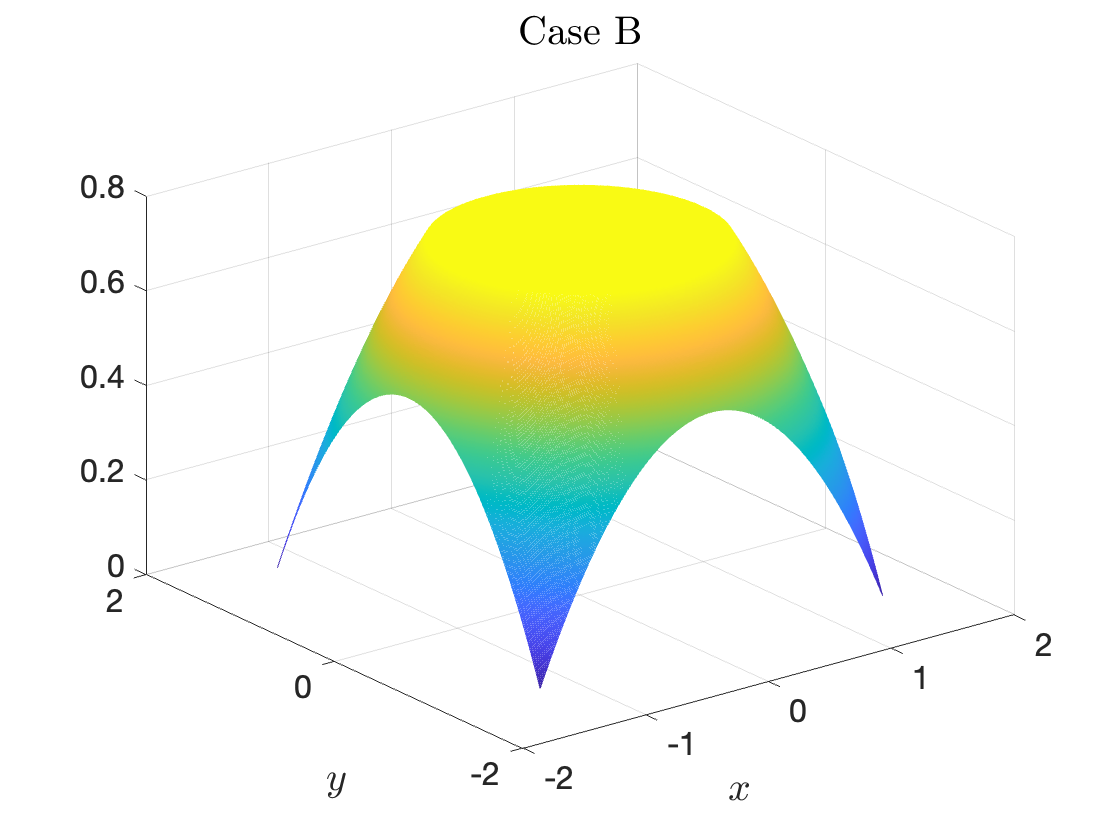}
\end{overpic}
\end{minipage}         
\begin{minipage}[b]
		{.32\textwidth}
		\centering
	\begin{overpic}[abs,width=\textwidth,unit=1mm,scale=.25]{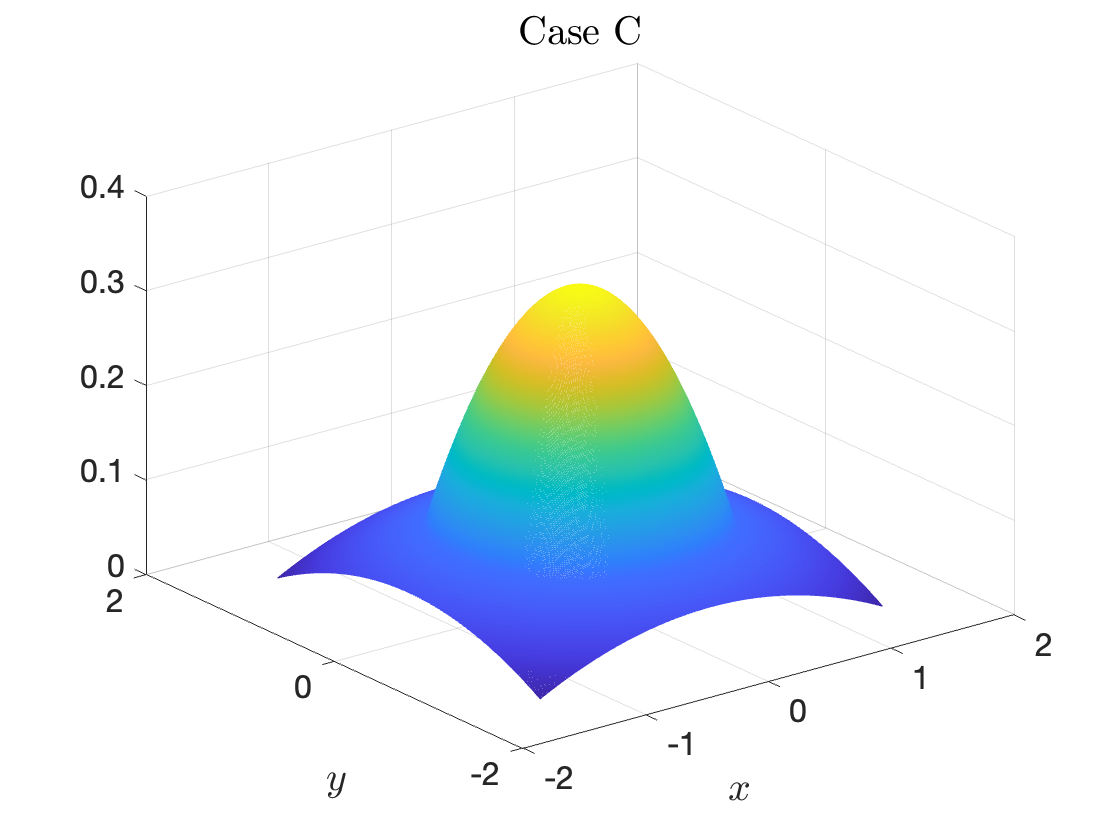}
\end{overpic}
\end{minipage}  
\begin{minipage}[b]
		{.32\textwidth}
		\centering
	\begin{overpic}[abs,width=\textwidth,unit=1mm,scale=.25]{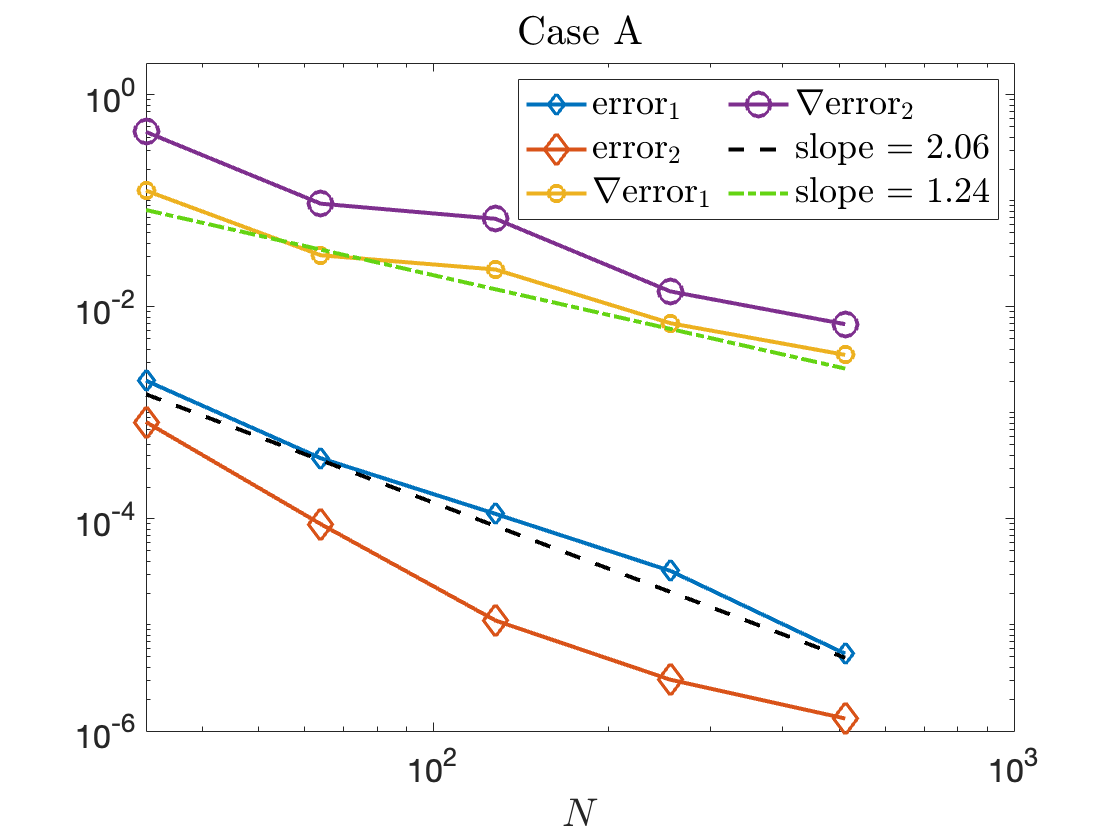}
\end{overpic}
\end{minipage}         
\begin{minipage}[b]
		{.32\textwidth}
		\centering
	\begin{overpic}[abs,width=\textwidth,unit=1mm,scale=.25]{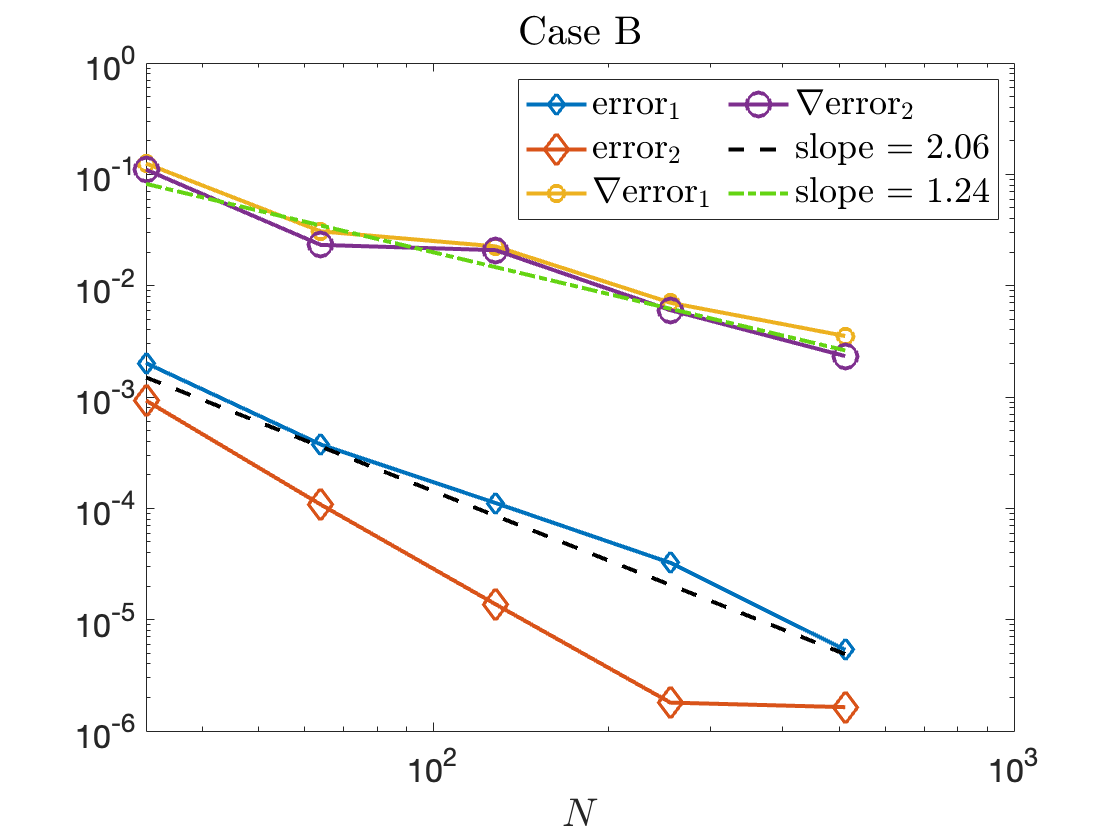}
\end{overpic}
\end{minipage}         
\begin{minipage}[b]
		{.32\textwidth}
		\centering
	\begin{overpic}[abs,width=\textwidth,unit=1mm,scale=.25]{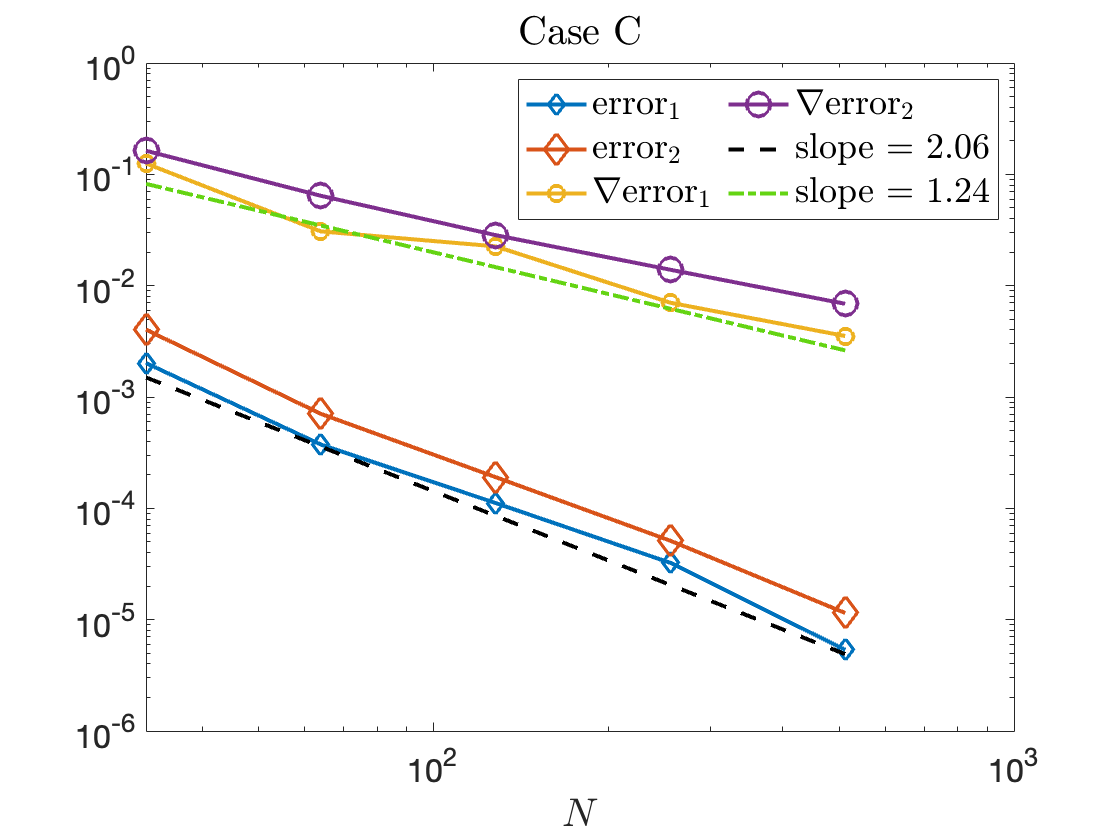}
\end{overpic}
\end{minipage}
    \caption{\textit{The first row shows the numerical solutions defined in Eq.~\eqref{eq:cases}, for the problem in Eq.~\eqref{eq:model}. The second row shows the convergence order of the respective numerical tests.}}
\label{fig:caseABC}
\end{figure}

Here, we perform numerical tests for a generic elliptic interface problem in two dimensions. Let $\Omega$ be a domain in $\mathbb R^2$ with external boundary $\partial \Omega_1$. We assume that $\Omega$ is subdivided into two subdomains $\Omega_i,\, i=1,2$ so that $\overline \Omega = \overline \Omega_1 \cup \overline \Omega_2$. The subdomains are separated by
a %Lipschitz 
continuous interface $\Gamma = \overline \Omega_1 \cap \overline \Omega_2$. Then, we
consider the following system:
\begin{subequations}
\label{eq:model}
    \begin{align}
        -\nabla \cdot \left( \beta_i \nabla u\right) & = f_i \qquad {\rm in }\, \Omega_i,\, i = 1,2 \\
        \left[\left|u\right|\right] & = g_D \qquad {\rm on }\, \Gamma \\
        \left[\left| \beta_i \nabla u \cdot \widehat n_i \right| \right] & = g_N \qquad {\rm on }\, \Gamma \\ 
        u & = g \qquad {\rm on} \, \partial \Omega_1
      \end{align}
\end{subequations}
where $f_i:\Omega_i \to \mathbb R, \, i = 1,2 $ and $\widehat n_i, \, i =1,2$ is the unit vector pointing out of $\Omega_i,\, i = 1,2$ and normal to $\Gamma$.
Our simulations involve various cases, including different combinations of coefficients between the two domains, and different shape domains. Specifically, we consider a square, a circle, and a flower, each
equipped with suitable Dirichlet boundary conditions. We first consider square $\Omega_1 = [-1.5,1.5]^2\setminus \Omega_2$ and $\Omega_2$ is defined by the following level set function 
\[ \phi(x,y) = \sqrt{x^2+y^2}-1.\]
In these cases, it is possible to find expressions for the exact solution of system~\eqref{eq:model} for $g_D = 0, g_N = 0$ and $g = u_1 $ on $\partial \Omega_1$. We test the following three cases:
\begin{subequations}
\label{eq:cases}
    \begin{align}
    {\rm Case\, A}& : \left\{
\begin{array}{l}
\displaystyle u^1_{\rm exa} = \frac{4 - x^2 - y^2}{4}, \quad \beta_1 = 1 \qquad {\rm in }\, \Omega_1  \\
\displaystyle u^2_{\rm exa} = \frac{31 - x^2 - y^2}{40}, \quad \beta_2 = 10  \qquad {\rm in }\, \Omega_2
\end{array}
\right. \\
 {\rm Case\, B}& : \left\{
\begin{array}{l}
\displaystyle u^1_{\rm exa} = \frac{4 - x^2 - y^2}{4}, \quad \beta_1 = 1 \qquad {\rm in }\, \Omega_1  \\
\displaystyle u^2_{\rm exa} = \frac{30001 - x^2 - y^2}{40000}, \quad \beta_2 = 10000  \qquad {\rm in }\, \Omega_2
\end{array}
\right. \\
 {\rm Case\, C}& : \left\{
\begin{array}{l}
\displaystyle u^1_{\rm exa} = \frac{4 - x^2 - y^2}{40}, \quad \beta_1 = 10 \qquad {\rm in }\, \Omega_1  \\
\displaystyle u^2_{\rm exa} = \frac{13 - 10x^2 - 10y^2}{40}, \quad \beta_2 = 1  \qquad {\rm in }\, \Omega_2
\end{array}
\right.
\end{align}
\end{subequations}
These examples are chosen from \cite{alshehri2024unfitted} to demonstrate the
robustness of our proposed method in dealing with different geometries, and the numerical results are shown in Fig.~\ref{fig:caseABC}. In the first row, we present the meshes of the numerical solutions for the three cases. In the second row, we show the error in the solution and its gradient, as defined in Eq.~\eqref{eq:error}. Since a quadrature rule is used to compute the integrals over the cells of the numerical solutions, quadrature errors must also be taken into account. As a result, the error curves exhibit a flattening behavior when reaching values on the order of $10^{-6}$. Table~\ref{tab:accuracy} reports the observed orders of convergence of the error for Case A (order$_1$ in $\Omega_1$ and order$_2$ in $\Omega_2$) as the parameter $\alpha$ in Algorithm~\ref{alg_snap} is varied. The fourth and fifth columns (i.e., $\alpha=(2,1.75)$ and $\alpha=(2,1.5)$) correspond to using different values of $\alpha$ in $\Omega_1$ and $\Omega_2$. These two tests show that the error decreases irregularly, achieving the second order only for larger values of $N$.

The next test explores the case in which the jump condition $g_N$ is nonzero. This scenario is included primarily for completeness, as the main applications of this work involve homogeneous interface conditions. We consider a flower shaped domain immersed in a square, where $\Omega_1 = [-1.5,1.5]^2 \setminus \Omega_2$ and $\Omega_2$ is defined by the following level set function 
\[ \phi(x,y) = r - 0.52 - \frac{Y^5 +5X^4Y-10X^2Y^3}{5r^5}, \quad r = \sqrt{X^2 + Y^2}, \, X = x-0.0\sqrt{3},\, Y = y-0.04\sqrt{2}.\]

\begin{figure}
    \centering
\begin{minipage}[b]
		{.49\textwidth}
		\centering
	\begin{overpic}[abs,width=\textwidth,unit=1mm,scale=.25]{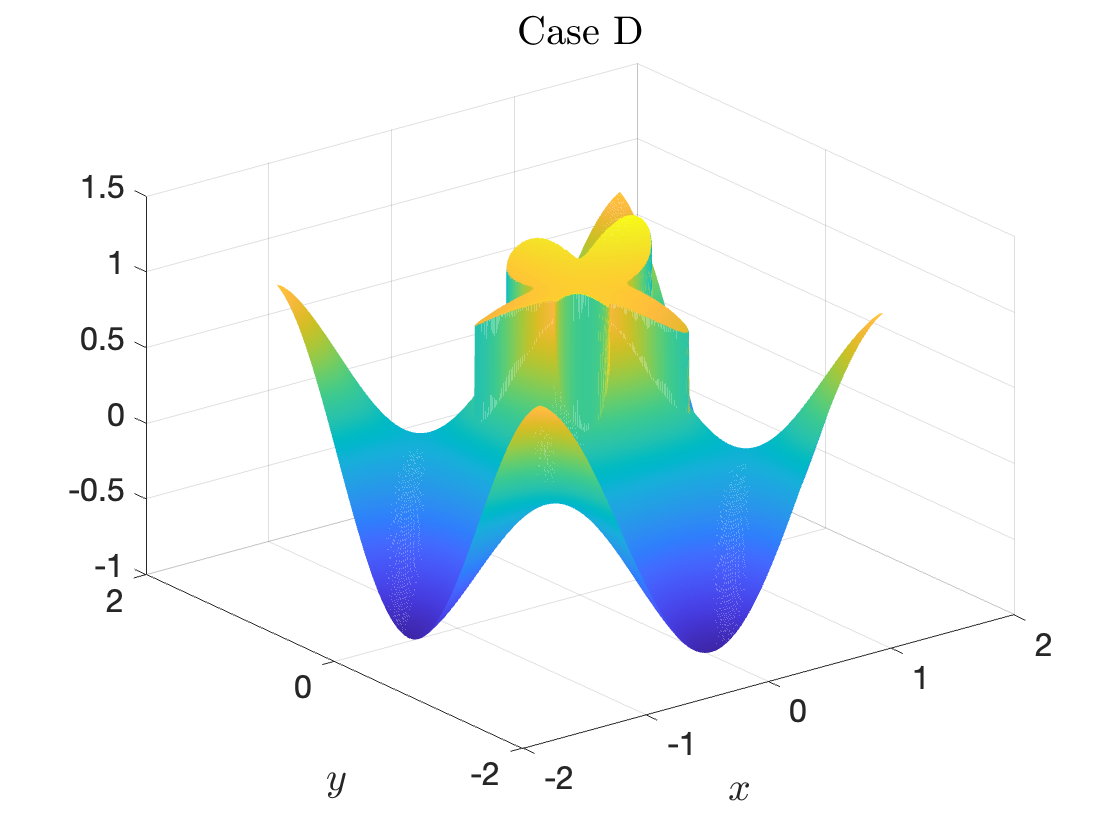}
    \put(2,65){(a)}
\end{overpic}
\end{minipage}         
\begin{minipage}[b]
		{.49\textwidth}
		\centering
	\begin{overpic}[abs,width=\textwidth,unit=1mm,scale=.25]{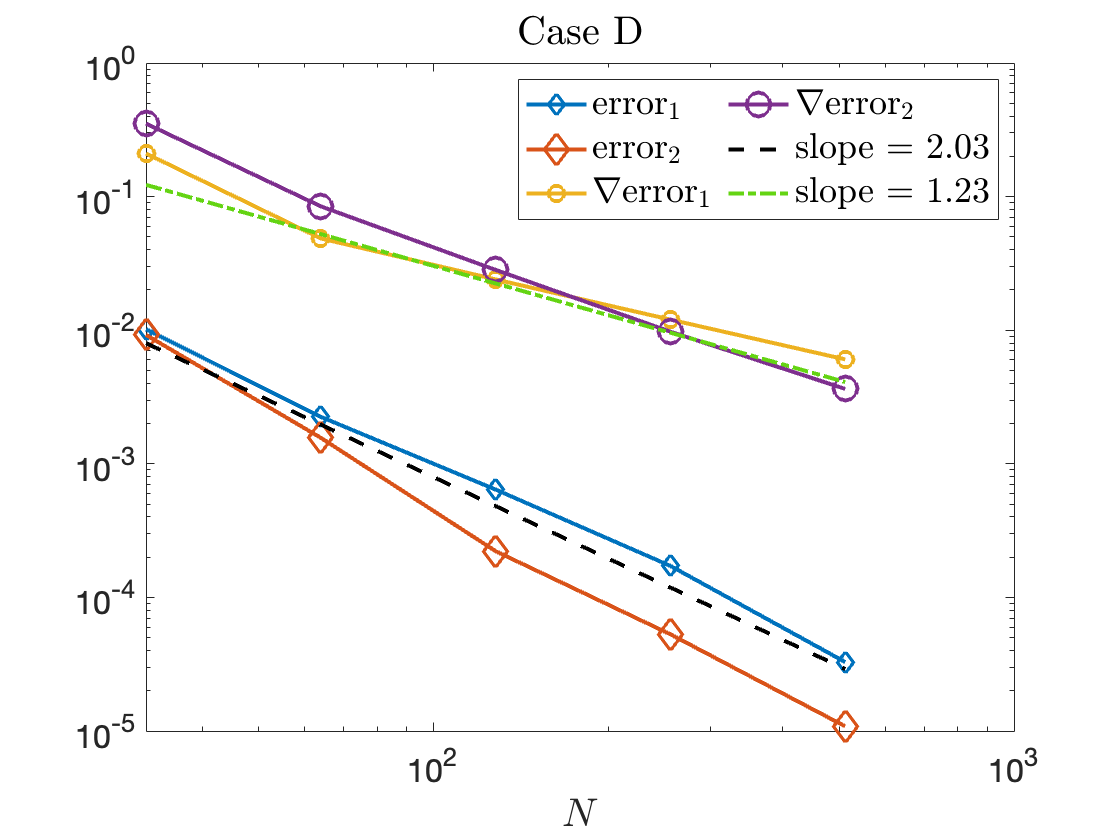}
    \put(2,65){(b)}
\end{overpic}
\end{minipage}          
    \caption{\textit{Panel (a) shows the numerical solution defined in Eq.~\eqref{eq:caseD}, for the problem in Eq.~\eqref{eq:model}. Panel (b) shows the convergence order of the respective numerical test.}}
\label{fig:caseD}
\end{figure}

The reference solutions we choose are the following 
\begin{subequations}
\label{eq:caseD}
    \begin{align}
    {\rm Case\, D}& : \left\{
\begin{array}{l}
\displaystyle u^1_{\rm exa} = \cos((2\pi x)/L)\cos((2\pi y)/L), \quad \beta_1 = 1 \qquad {\rm in }\, \Omega_1  \\
\displaystyle u^2_{\rm exa} = \exp(xy), \quad \beta_2 = 2  \qquad {\rm in }\, \Omega_2
\end{array}
\right.
\end{align}
\end{subequations}
and the numerical results are shown in Fig.~\ref{fig:caseD}. Panel (a) shows the meshes of the numerical solution. Panel (b) shows the error in the solution and its gradient, as defined in Eq.~\eqref{eq:error}.

\begin{table}[H]
\centering
\caption{\textit{Order of accuracy of Case A, changing the parameter $\alpha$ in Algorithm~\ref{alg_snap}.}}
\label{tab:accuracy}
\begin{tabular}{||c||c|c||c|c||c|c||c|c||c|c||}
\hline 
& \multicolumn{2}{c||}{$\alpha = 2$}
& \multicolumn{2}{c||}{$\alpha = 1.75$}
& \multicolumn{2}{c||}{$\alpha = 1.5$} 
& \multicolumn{2}{c||}{$\alpha = (2,1.75)$} 
& \multicolumn{2}{c||}{$\alpha = (2,1.5)$} \\ \hline
$N$ & order$_1$ & order$_2$
 & order$_1$ & order$_2$
& order$_1$ & order$_2$
 & order$_1$ & order$_2$
 & order$_1$ & order$_2$\\
\hline
32 & 2.02 & 5.20 & 0.93 & 2.63 & 0.94 & 2.19 & 0.66 & 1.13 & 2.51 & 3.16 \\ \hline
64 & 2.05 & 1.46 & 2.38 & 5.86 & 0.86 & 2.54 & 5.10 & 5.40 & 2.53 & 1.08 \\ \hline
128 & 2.03 & 2.34 & 3.48 & 1.84 & 0.53 & 0.0 & 1.28 & 1.88 & 2.08 & 0.013 \\ \hline
256 & 1.99 & 1.89 & 2.16 & 1.70 & 2.10 & 3.13 & 1.41 & 1.72 & 1.40 & 6.28 \\ \hline
512 & 1.97 & 2.04 & 1.56 & 1.92 & 1.07 & 1.38 & 3.96 & 2.13 & 3.50 & 1.18 \\ \hline
1024 & 1.99 & 2.16 & 2.04 & 1.93 & 1.90 & 2.10 & 1.48 & 1.96 & 1.92 & 1.80 \\ 
\hline
\end{tabular}
\end{table}

\section{Two interfaces problem and Numerical results}
\label{SEC:num_res}
We now consider the problem introduced in Section \ref{SEC:math_model}, in which the electrostatic potential is affected by the presence of two interfaces. The numerical simulations are carried out by solving the drift-diffusion equations (\ref{eq:currents}-\ref{eq:currentsnp}) self-consistently coupled with the Poisson equation \eqref{eq:PoissonPhi}. The resulting nonlinear system is solved by an iterative procedure until self-consistency between the carrier densities and the electrostatic potential is achieved.

%We choose an initial guess for the electron and hole densities, $n^0$ and $p^0$ respectively, and 
%we use the intrinsic densities, namely
% \begin{subequations}
%     \begin{align}
%         n^0 &= n(0), \\
%         p^0 &= p(0),
%     \end{align}
% \end{subequations}
% where \(n(0)\) and \(p(0)\) are obtained from the relations \eqref{eq:np_relations}.
%the corresponding electrostatic potential \(\varphi^0\) is then computed by solving the Poisson equation \eqref{eq:PoissonPhi}. 

At a generic iteration \(k\), given the densities \(n^k\) and \(p^k\), the Poisson equation \eqref{eq:PoissonPhi} is solved to obtain the electrostatic potential \(\varphi^k\). The drift-diffusion equations (\ref{eq:currents}-\ref{eq:currentsnp}) are then solved as described in Section \ref{sec:numerical}, yielding the quasi-Fermi energy profiles \(\left(\varepsilon^{(n)}\right)^k\) and \(\left(\varepsilon^{(p)}\right)^k\). New carrier densities are subsequently computed as
\begin{subequations}
    \begin{align}
        n^{\mathrm{new},k} &= n\left(\left(\varepsilon^{(n)}\right)^k + e\varphi^k\right), \\
        p^{\mathrm{new},k} &= p\left(\left(\varepsilon^{(p)}\right)^k + e\varphi^k\right).
    \end{align}
\end{subequations}
To improve the stability of the self-consistent iteration, we use a fixed under-relaxation step. Starting with a couple of initial guesses $(n^0,p^0)$, the updated densities are defined by
\begin{subequations}
    \begin{align}
        n^{k+1} &= \xi\, n^k + (1-\xi)n^{\mathrm{new},k}, \\
        p^{k+1} &= \xi\, p^k + (1-\xi)p^{\mathrm{new},k},
    \end{align}
\end{subequations}
with \(0<\xi<1\), until the difference between two successive approximations falls below a specific tolerance {\tt toll}. This corresponds to a damped fixed-point iteration with relaxation parameter \(1-\xi\) \cite{DederichsZeller1983}.

The computational domain consists of a three-layer sandwich structure composed of two inactive regions, denoted by $\Omega_1$ and $\Omega_2$, separated by an active graphene layer, denoted by $\Omega_{\mathrm{gr}}$. In this study, we investigate the effect of varying the thickness of the graphene layer, defined as
\[
t_{\mathrm{gr}} = y_3 - y_2.
\]
To conduct this comparison, we consider three values of $t_{\mathrm{gr}}$. As a consequence, the total height $H$ of the computational domain varies across the numerical experiments.

Each layer is represented through a level set function describing a rectangular region centered at $(x_c,y_c)$. Let
\[
\widetilde{x}=x-x_c, \qquad \widetilde{y}=y-y_c,
\]
and define
\begin{equation}
\label{eq:levelset1}
\phi_i(x,y)
=
\max\!\left(|\widetilde{x}|-r_x^{(i)},\,|\widetilde{y}|-r_y^{(i)}\right),
\qquad i=1,\mathrm{gr},2,
\end{equation}
where $r_x^{(i)}$ and $r_y^{(i)}$ denote the half-width and half-height of the corresponding layer.

The three subdomains are then given by
\begin{equation}
\Omega_i
=
\left\{
(x,y)\in\mathbb{R}^2 :
\phi_i(x,y)<0
\right\},
\qquad
i=1,\mathrm{gr},2.
\end{equation}

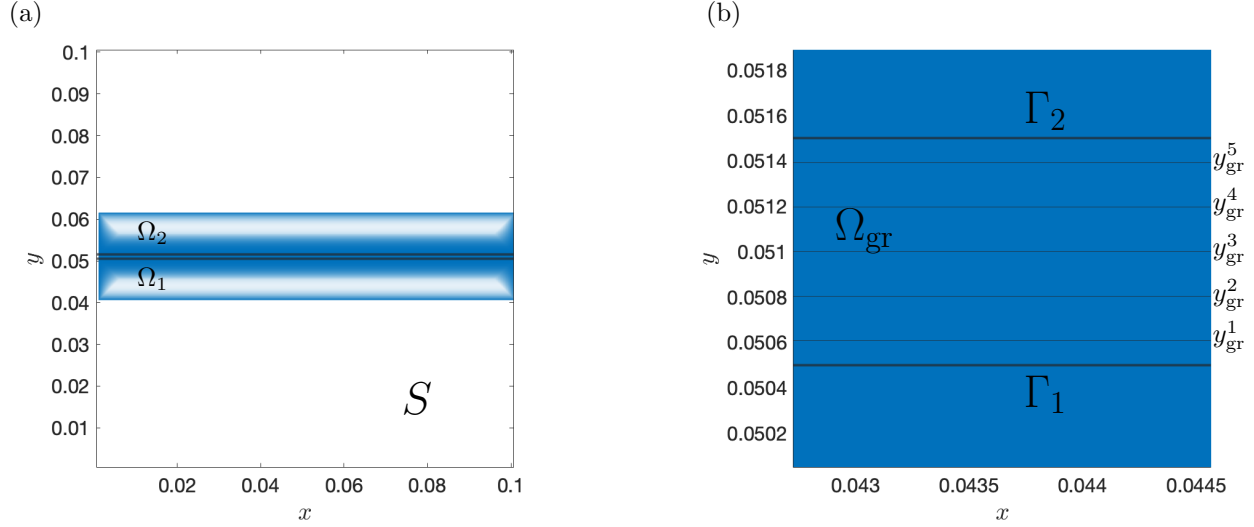
\begin{figure}[ht]
    \centering
\begin{minipage}[b]
		{.49\textwidth}
		\centering
	\begin{overpic}[abs,width=\textwidth,unit=1mm,scale=.25]{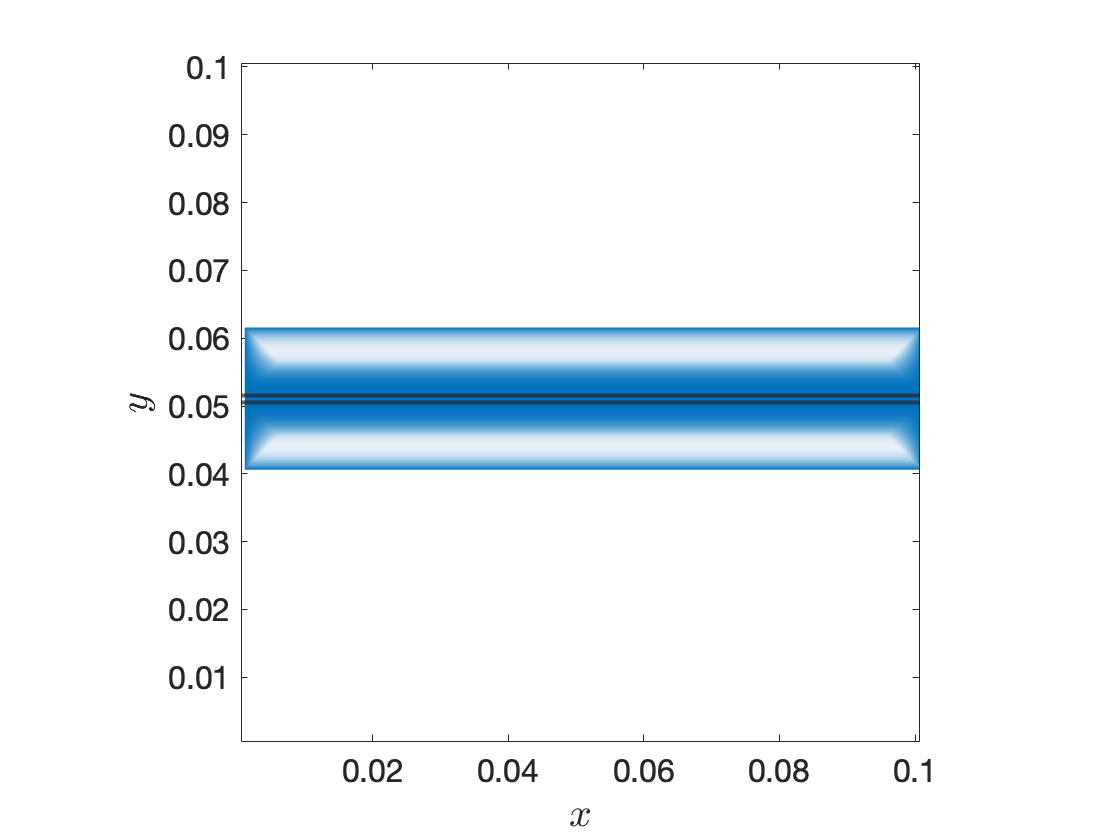}
    \put(8,67){(a)}
    \put(60,15){\LARGE $S$}
    \put(25,32){$\Omega_1$}
    \put(25,38){$\Omega_2$}
\end{overpic}
\end{minipage}         
\begin{minipage}[b]
		{.49\textwidth}
		\centering
	\begin{overpic}[abs,width=\textwidth,unit=1mm,scale=.25]{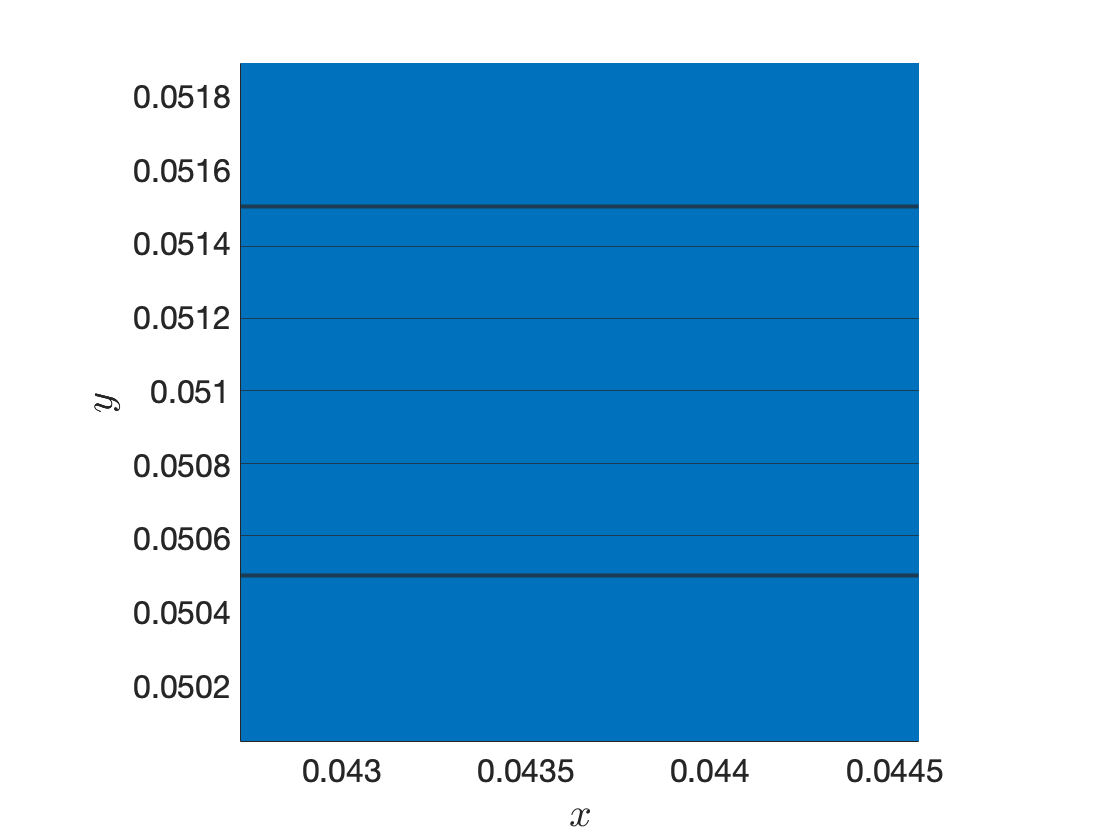}
   \put(8,67){(b)}
   \put(25,38){\LARGE  $\Omega_{\rm gr}$}
    \put(50,16){\LARGE  $\Gamma_1$}
    \put(50,53.5){\LARGE  $\Gamma_2$}
    \put(75,36){$y^3_{\rm gr}$}
    \put(75,24.5){$y^1_{\rm gr}$}
    \put(75,30){$y^2_{\rm gr}$}
    \put(75,42){$y^4_{\rm gr}$}
    \put(75,48){$y^5_{\rm gr}$}
\end{overpic}
\end{minipage}          
    \caption{\textit{(a) Mesh of the region where the product of the level-set functions $\phi_1, \phi_2, \phi_{\rm gr}$ (defined in Eq.~\eqref{eq:levelset1}) is negative, which corresponds to the computational domain. Zoom of the central layer $\Omega_{\rm gr}$ in panel (b). Parameters: $t_{\rm gr} = 1$nm, $N = 520$.}}
    \label{fig:levelset}
\end{figure}

All layers are rectangular and aligned along the same centerline, differing only in their vertical position and thickness, see Fig.~\ref{fig:levelset}. The computational domain is symmetric with respect to the horizontal centerline passing through the center of the graphene layer $\Omega_{\mathrm{gr}}$, whose position is kept fixed throughout the analysis. Consequently, any modification of the central layer is performed in a way that preserves this symmetry with respect to the midline of $\Omega_{\mathrm{gr}}$. The parameters of the numerical simulations are shown in Table~\ref{table_parameters}. 

\begin{figure}[h]
    \centering
\begin{minipage}[b]
		{.49\textwidth}
		\centering
	\begin{overpic}[abs,width=\textwidth,unit=1mm,scale=.25]{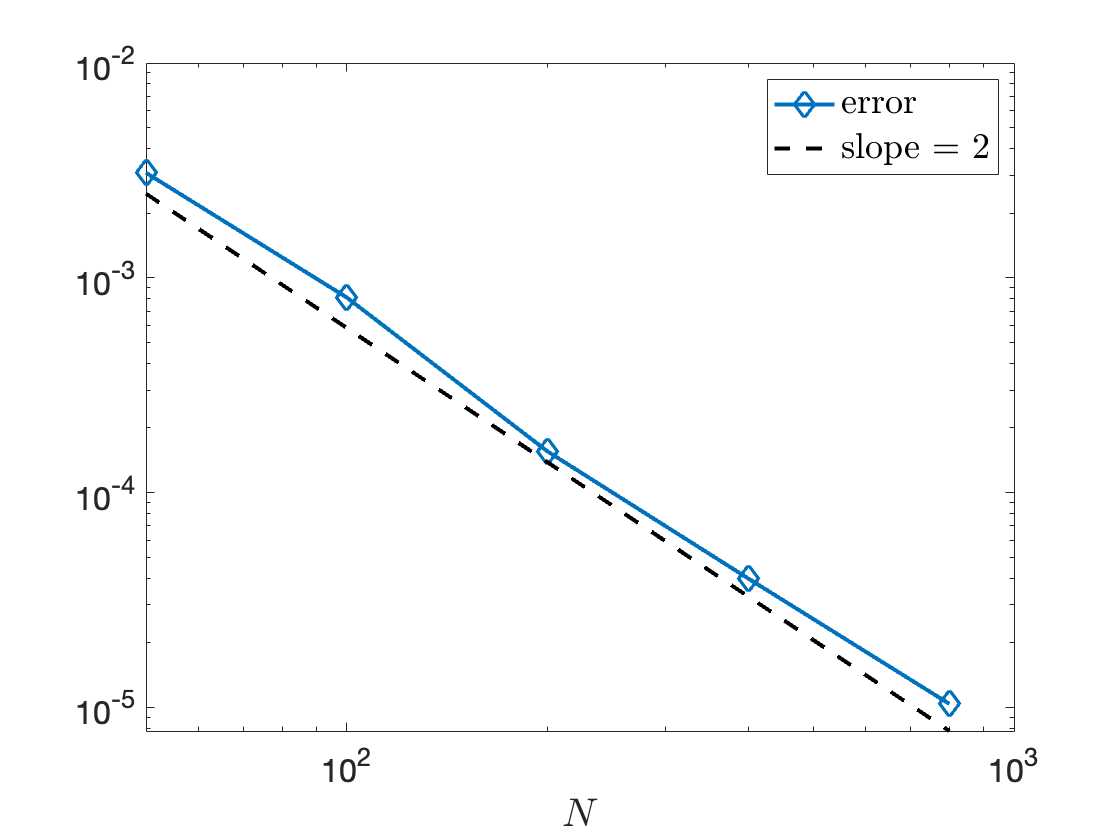}
    \put(2,65){(a)}
\end{overpic}
\end{minipage} 
\begin{minipage}[b]
		{.49\textwidth}
		\centering
	\begin{overpic}[abs,width=\textwidth,unit=1mm,scale=.25]{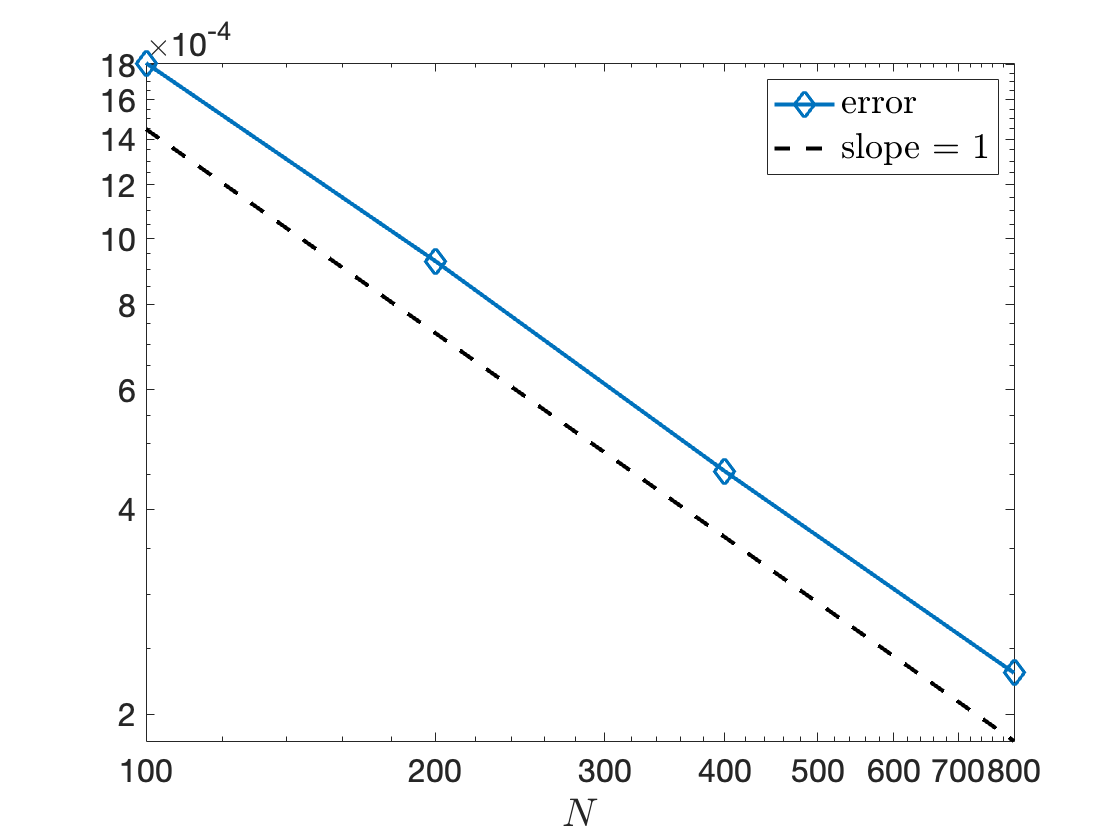}
    \put(2,65){(b)}
\end{overpic}
\end{minipage} 
    \caption{\textit{Convergence order of the numerical accuracy of the potential $\varphi$ in the three-layer domain defined by Eq.~\eqref{eq:levelset1} after one iteration. In panel (a), only the two-dimensional  elliptic equation is solved in $\Omega$. In panel (b), the full coupled system is considered, combining the elliptic problem with the one-dimensional drift-diffusion equation. In this test, the order is calculated with Richardson extrapolation technique.}}
    \label{fig:accuracy_3layers}
\end{figure}

Fig.~\ref{fig:accuracy_3layers} presents the convergence study for the numerical solution of the potential $\varphi$, in the case of a three layers domain $\Omega$ with two interfaces $\Gamma_1$ and $\Gamma_2$. The results are shown after one iteration, and are obtained by solving only the elliptic problem in panel (a), while panel (b) displays the results for the fully coupled drift-diffusion-Poisson system. The obtained convergence rates confirm the expected order of accuracy, as the nonlinear terms arising from the transport model reduce the overall convergence order to one in panel (b). In the absence of an exact solution, we apply the Richardson extrapolation technique (see, e.g., \cite{richardson1911ix}) to estimate the order of the method. 


\begin{figure}
    \centering
\begin{minipage}[b]
		{.49\textwidth}
		\centering
	\begin{overpic}[abs,width=\textwidth,unit=1mm,scale=.25]{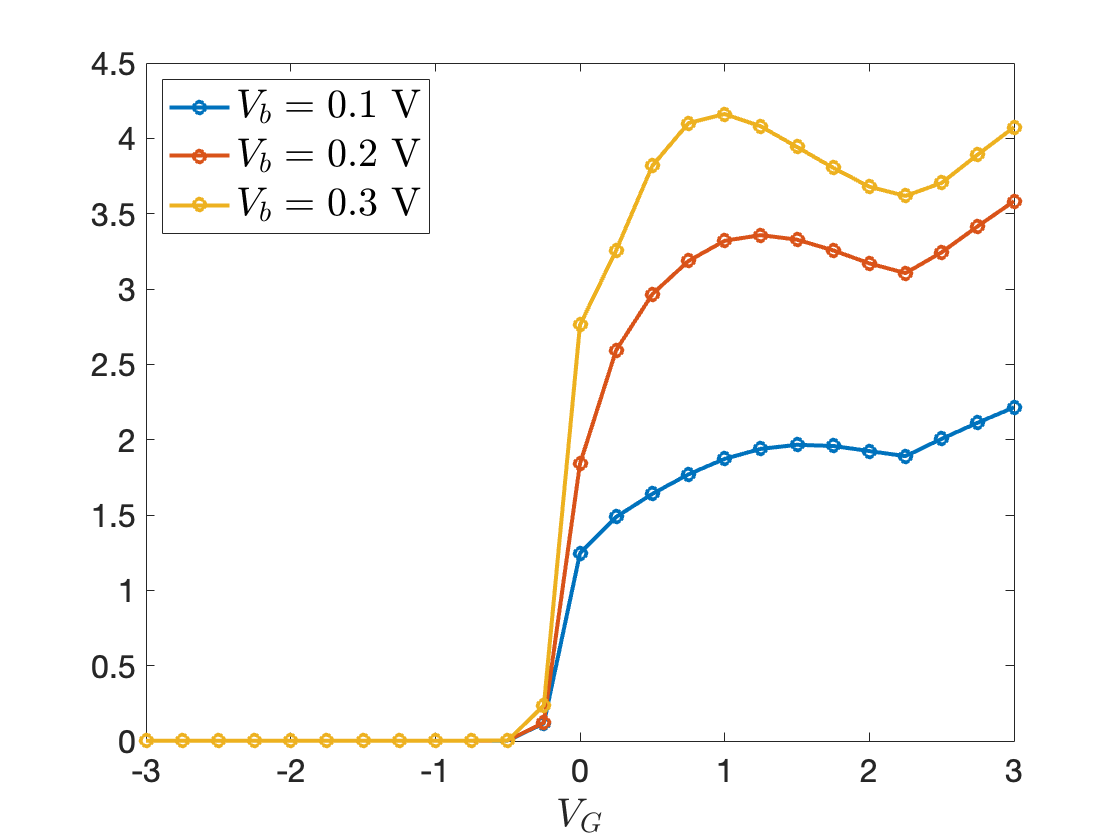}
    \put(2,65){(a)}
\end{overpic}
\end{minipage}         
\begin{minipage}[b]
		{.49\textwidth}
		\centering
	\begin{overpic}[abs,width=\textwidth,unit=1mm,scale=.25]{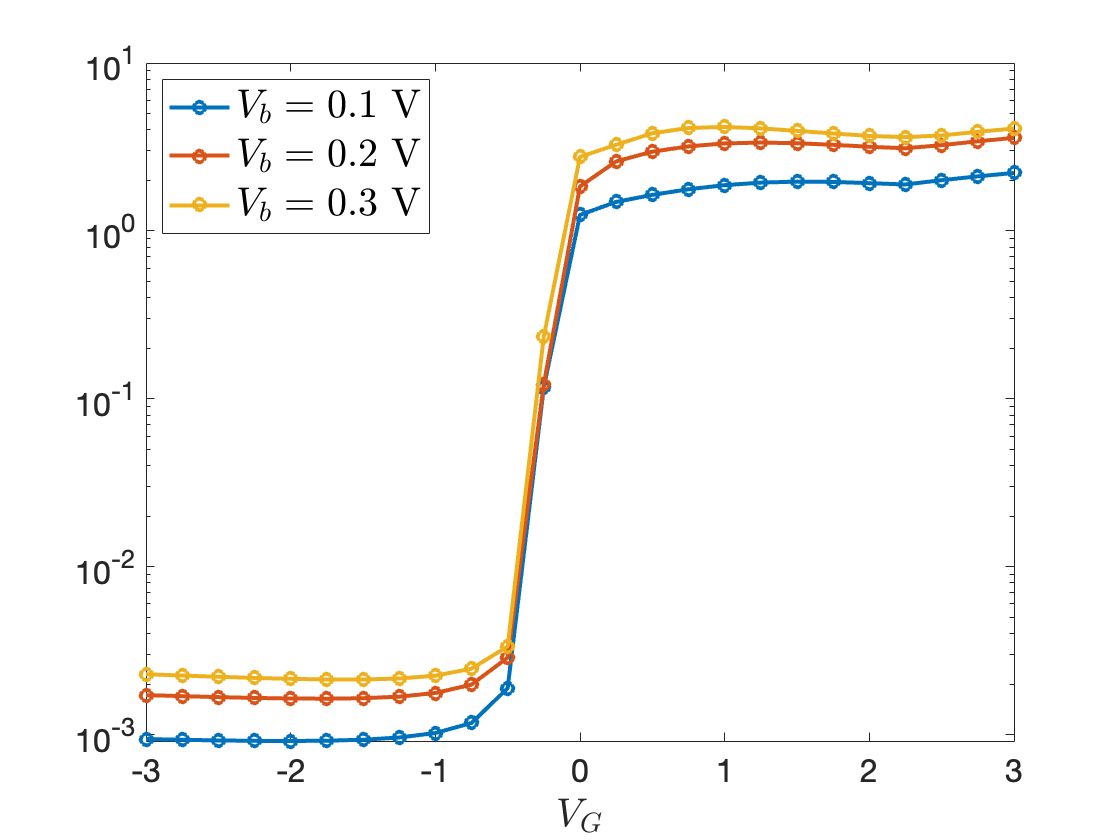}
    \put(2,65){(b)}
\end{overpic}
\end{minipage}          
    \caption{\textit{Characteristic curves of the simulated GFET in linear scale in panel (a) and logarithmic one in panel (b). Parameters: $t_{\mathrm{gr}} = 1$nm, initial values $n^0 = p^0 = 10$, number of points $N = 520$, number of lines in the graphene layer $N_{\rm gr} = 5$, snapping parameter $\alpha = 2$, $\xi= 10^{-4}$ and {\tt toll} = $10^{-4}$.}}
    \label{fig:charact}
\end{figure}

Fig.~\ref{fig:charact} shows the transfer characteristics of the simulated GFET, defined in Eq.~\eqref{eq:Jav}, for three values of the drain-source bias voltage, namely $V_b=0.1$, $0.2$, and $0.3$ V. The same data are reported both in linear and logarithmic scale in order to highlight, respectively, the OFF-state and ON-state behavior. For negative gate voltages, the current remains strongly suppressed, while a sharp increase is observed when $V_G$ approaches the threshold region. This indicates that the electrostatic control exerted by the gates is effective in driving the device from a low-conducting state to an ON state. Such behavior is particularly relevant in monolayer graphene devices, where the absence of an energy gap usually leads to ambipolar conduction and makes the definition of a clear OFF region difficult.

As expected, the ON-state current increases with the drain-source bias voltage, due to the larger longitudinal driving field along the graphene channel. In the positive-gate region the current reaches a quasi-saturated regime, with only moderate variations as $V_G$ is further increased. This behavior can be ascribed to the combined effect of gate-induced carrier-density modulation, field-dependent mobility, and the self-consistent electrostatic potential obtained from the coupled transport-Poisson model. The logarithmic plot shows that the ratio between the ON current and the residual OFF current is of several orders of magnitude, confirming that the adopted device configuration mitigates the typical ambipolar response of gapless graphene and yields transfer curves with a conventional FET-like behavior.

A weak non-monotonicity is observed in the ON region of the transfer characteristics. Non-monotone or distorted transfer characteristics are not unusual in low-dimensional field-effect devices, since the gate voltage can affect not only the channel charge density but also the contact regions, the local electrostatic profile, and the effective transport parameters. This has been observed, for example, in graphene FETs with contact-induced doping or gate-dependent contact resistance, and is consistent with the general view that 2D/metal and 2D/dielectric interfaces play a decisive role in the electrostatics and transport of atomically thin transistors \cite{Nouchi2008,DiBartolomeo2013,Akinwande2019}.

%Fig.~\ref{fig:elpot} shows a typical shape of the 2D electrostatic potential. 
%To analyze the behavior of the solution within the active region, we examine the horizontal profiles of the electrostatic potential along the grid lines contained in the graphene layer,
%\[
%\varphi(x,y_{\rm gr}^j), \qquad j=1,\ldots,N_{\rm gr}.
%\]
%Fig.~\ref{fig:xprofile} shows these profiles for the discrete locations \(y_{\rm gr}^j\) inside \(\Omega_{\rm gr}\). Since the graphene layer is extremely thin, the comparison of the curves provides insight into the variation of the potential across the active region and allows us to assess the influence of the gate voltage on the in-plane potential distribution.

\begin{figure}
    \centering
\begin{minipage}[b]
		{.49\textwidth}
		\centering
	\begin{overpic}[abs,width=\textwidth,unit=1mm,scale=.25]{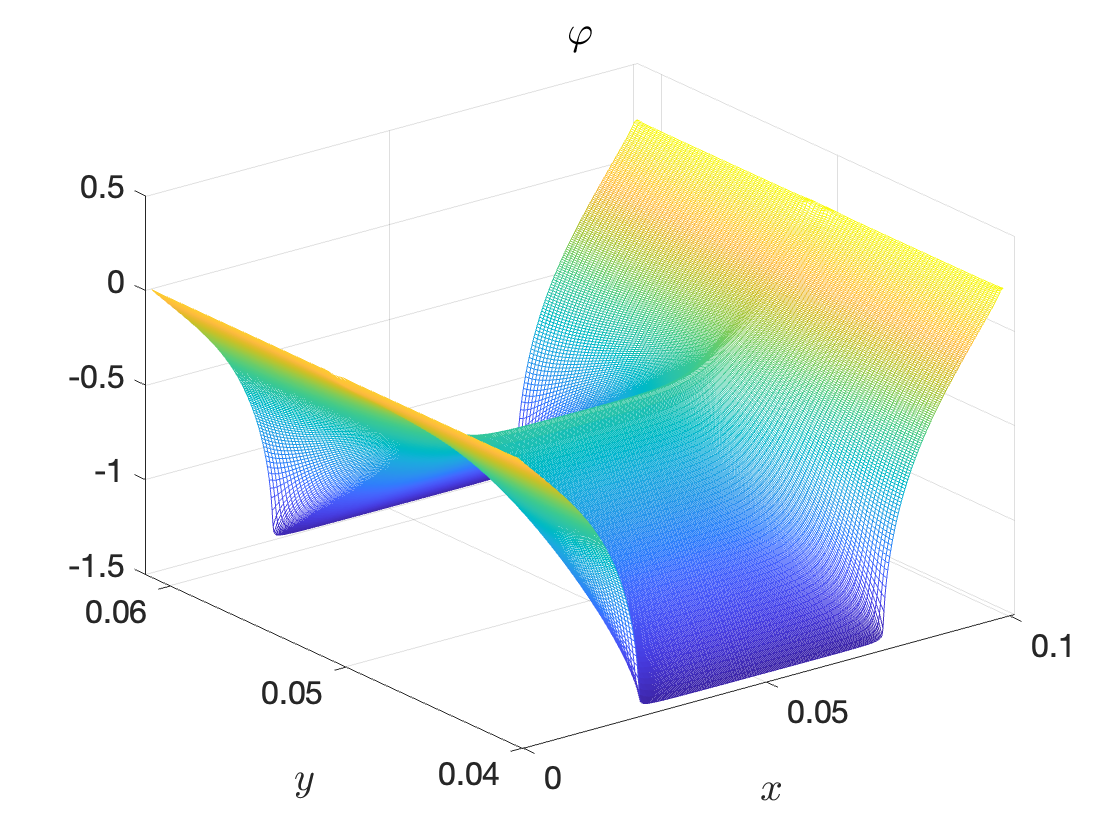}
    \put(2,65){(a)}
\end{overpic}
\end{minipage}         
\begin{minipage}[b]
		{.49\textwidth}
		\centering
	\begin{overpic}[abs,width=\textwidth,unit=1mm,scale=.25]{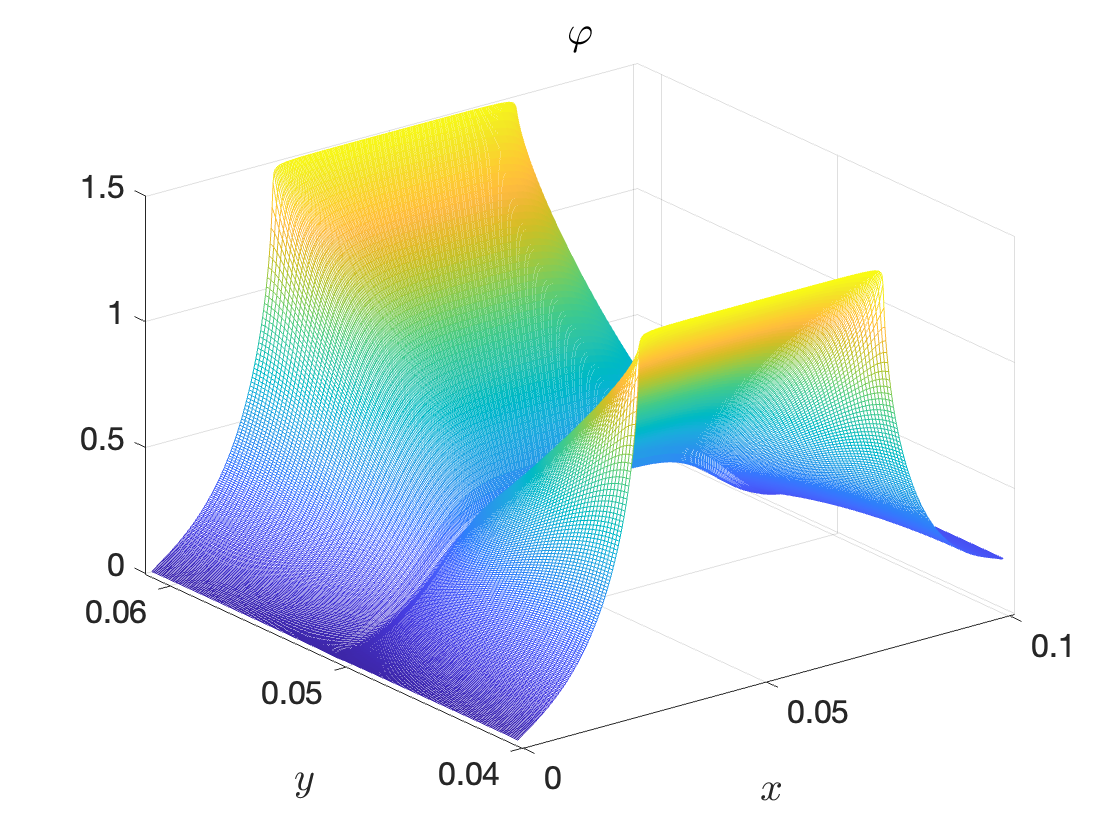}
    \put(2,65){(b)}
\end{overpic}
\end{minipage}          
    \caption{\textit{Electrostatic potential $\varphi$ of the simulated GFET. Parameters: $t_{\mathrm{gr}} = 1$nm, initial values $n^0 = p^0 = 10$, number of points $N= 520$, number of lines in the graphene layer $N_{\rm gr} = 5$, snapping parameter $\alpha = 2$, $\xi= 10^{-4}$, {\tt toll} = $10^{-4}$, $V_b = 0.2 V$ gate voltage $V_G = -1.5 V$ in panel (a) and $V_G = 1.5 V$ in panel (b).}}
    \label{fig:elpot}
\end{figure}

Fig.~\ref{fig:elpot} displays the two-dimensional electrostatic potential for $V_b=0.2\,\mathrm{V}$ and different values of the gate voltage: in panel (a) $V_G = -1.5$V and in panel (b) $V_G = 1.5$V. The potential is strongly affected by the gate bias: a negative $V_G$ lowers the potential in the central part of the channel, while a positive $V_G$ raises it and produces a broad plateau in the gated region. This behavior reflects the combined action of the source/drain boundary conditions and of the upper and lower gates acting through the two oxide/graphene interfaces.

\begin{figure}
    \centering
\begin{minipage}[b]
		{.49\textwidth}
		\centering
	\begin{overpic}[abs,width=\textwidth,unit=1mm,scale=.25]{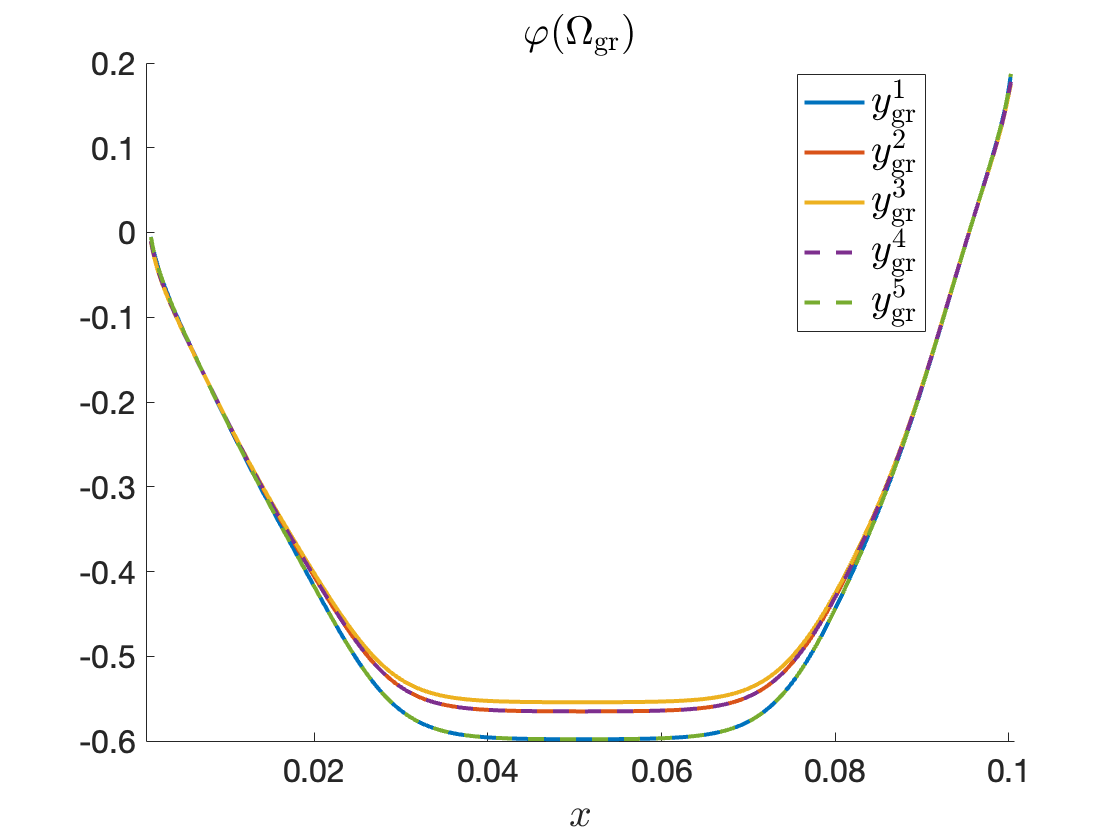}
    \put(2,65){(a)}
\end{overpic}
\end{minipage}         
\begin{minipage}[b]
		{.49\textwidth}
		\centering
	\begin{overpic}[abs,width=\textwidth,unit=1mm,scale=.25]{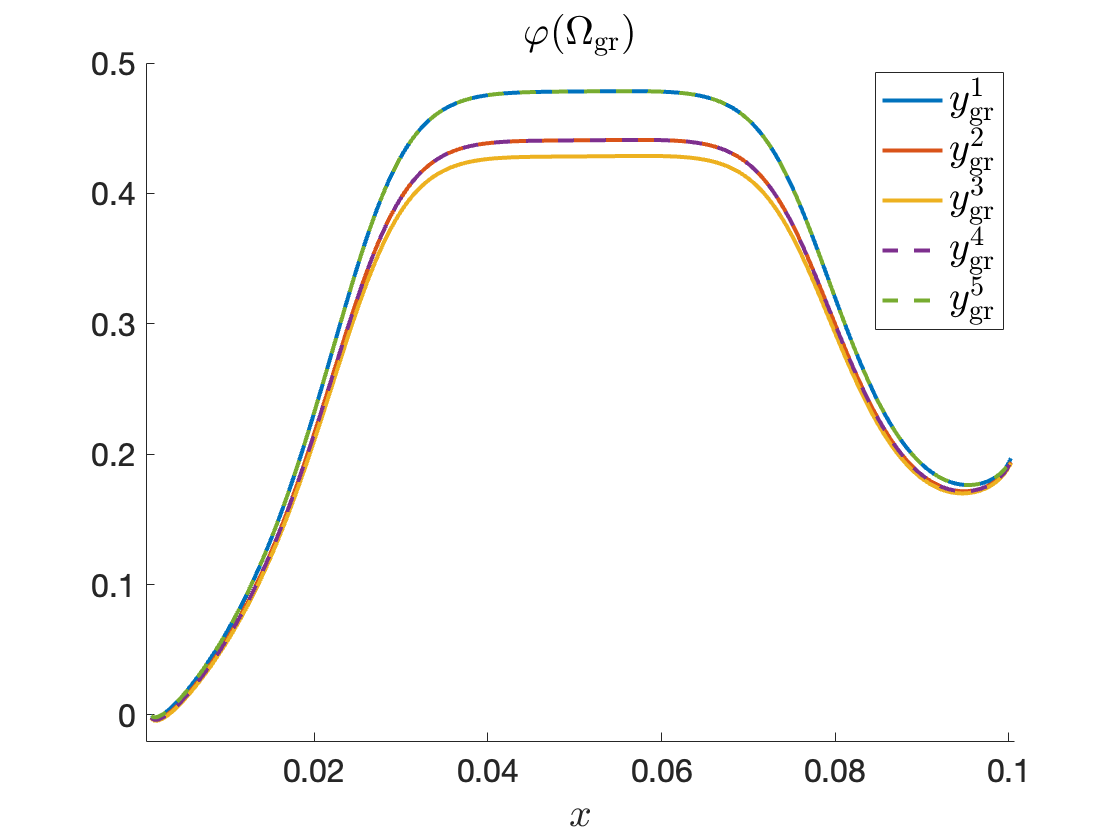}
    \put(2,65){(b)}
\end{overpic}
\end{minipage}          
    \caption{\textit{Profiles in the $x$-direction of the numerical solution $\varphi(x, y = y^j_{\rm gr}), \, j = 1,\cdots, 5$, in the graphene layer $\Omega_{\rm gr}$. 
    Parameters: $t_{\mathrm{gr}} = 1$nm, initial values $n^0 = p^0 = 10$, number of points $N = 520$, number of lines in the graphene layer $N_{\rm gr} = 5$, snapping parameter $\alpha = 2$, $\xi= 10^{-4}$, {\tt toll} = $10^{-4}$, $V_b = 0.2 V$ gate voltage $V_G = -1.5 V$ in panel (a) and $V_G = 1.5 V$ in panel (b).}}
    \label{fig:xprofile}
\end{figure}

The corresponding horizontal profiles inside the graphene layer are shown in Fig.~\ref{fig:xprofile}, where
\[
    \varphi(x,y_{\rm gr}^j), \qquad j=1,\ldots,N_{\rm gr},
\]
is plotted for the different grid lines belonging to $\Omega_{\rm gr}$. The curves are nearly superposed, as expected from the very small thickness of the graphene layer, but their small separation indicates that the transverse variation of the electrostatic potential is not completely negligible. This confirms the relevance of solving the full two-dimensional Poisson problem across the oxide/graphene/oxide structure. In particular, the gate voltage modifies not only the potential level in the channel, but also the longitudinal profile and hence the electric field, which is the quantity directly entering the drift-diffusion current.

%To further investigate the effect of the discretization of the graphene layer, we analyze the vertical profile of the electrostatic potential along the line passing through the center of the device, namely
%\[
%\varphi\left(\frac{x_1+x_4}{2},y\right).
%\]
%Fig.~\ref{fig:yprofile} reports this quantity for three different thicknesses \(t_{\rm gr}\) of the graphene layer \(\Omega_{\mathrm{gr}}\). The total number of grid points is kept fixed (\(N=820\)) in all simulations. Therefore, as \(t_{\rm gr}\) decreases, the number of grid points \(N_{\rm gr}\) contained in the graphene layer is progressively reduced. This configuration allows us to assess the sensitivity of the numerical solution to the resolution of the active layer while maintaining the same overall computational cost. 

\begin{figure}
    \centering
\begin{minipage}[b]
		{.49\textwidth}
		\centering
	\begin{overpic}[abs,width=\textwidth,unit=1mm,scale=.25]{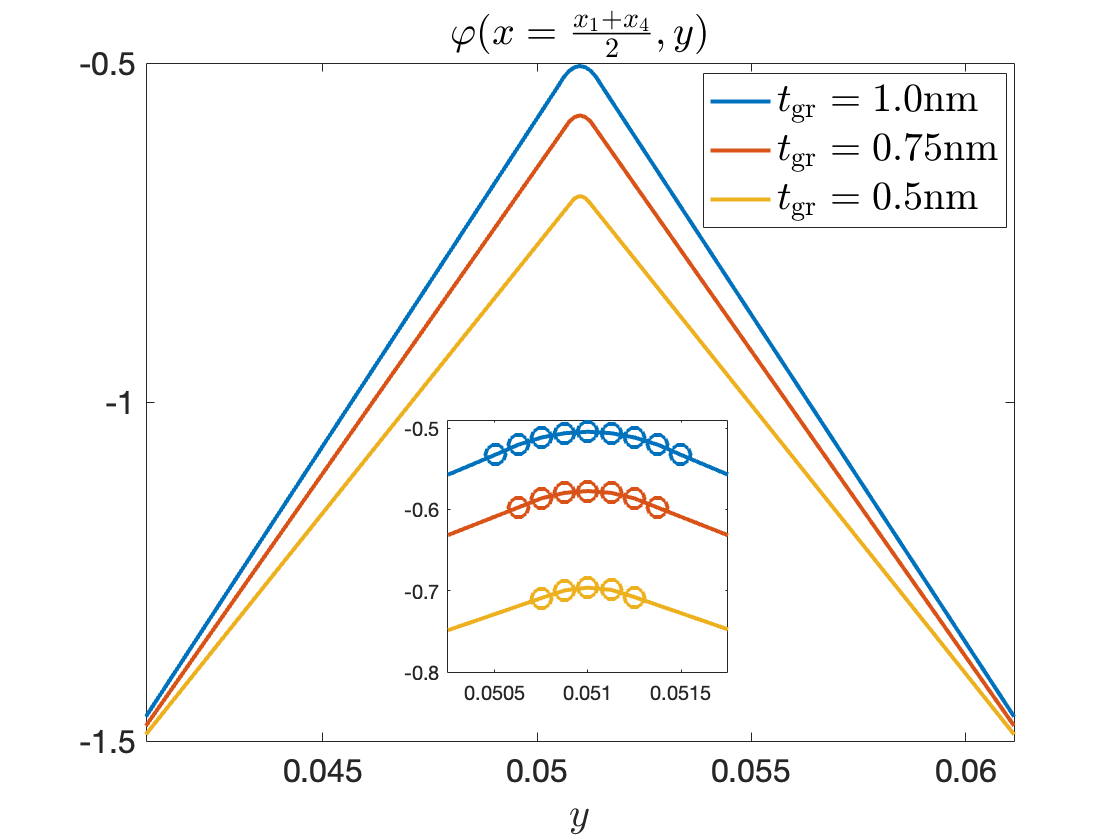}
    \put(2,65){(a)}
\end{overpic}
\end{minipage}         
\begin{minipage}[b]
		{.49\textwidth}
		\centering
	\begin{overpic}[abs,width=\textwidth,unit=1mm,scale=.25]{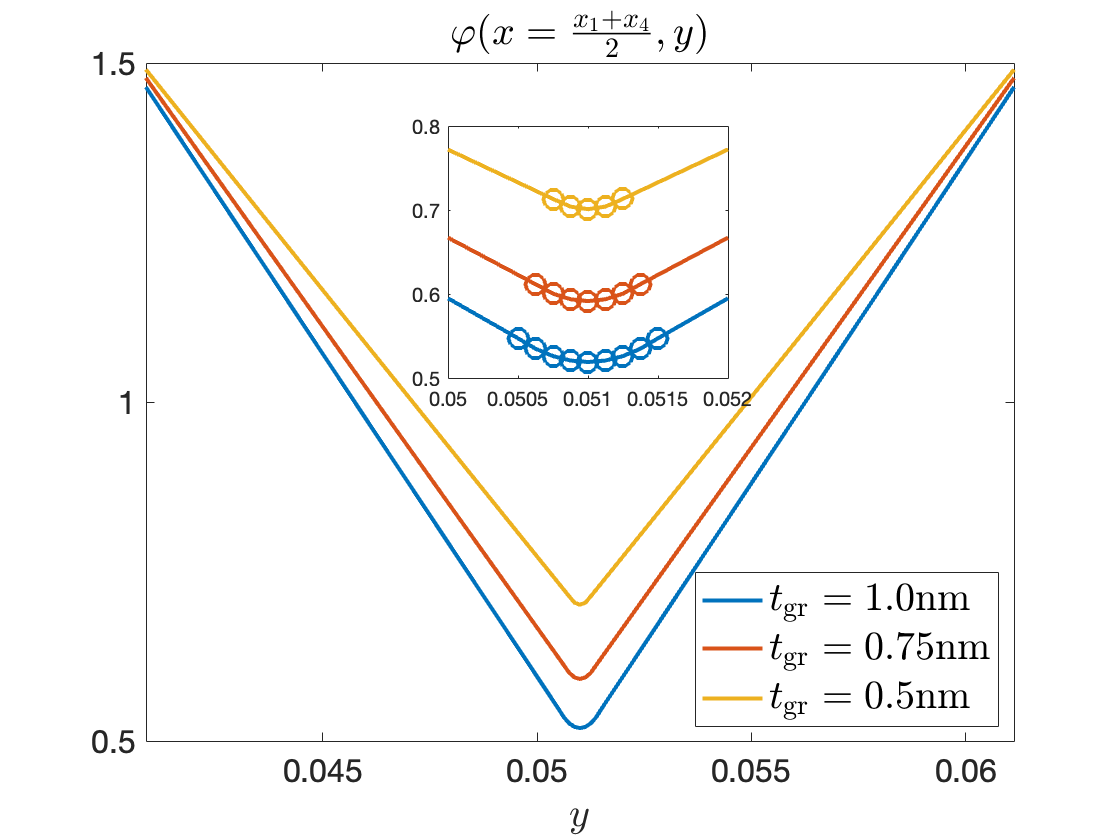}
     \put(2,65){(b)}
\end{overpic}
\end{minipage}          
   \caption{\textit{Profiles in the $y$-direction of the numerical solution $\varphi(x = \frac{x_1+x_4}{2}, y)$ for three configurations corresponding to different thicknesses $t_{\rm gr}$ of the graphene layer $\Omega_{\mathrm{gr}}$. The total number of grid points, $N = 820$, is fixed across all cases. As the thickness of $\Omega_{\mathrm{gr}}$ decreases, the number of points $y_{\mathrm{gr}}^{j}$, $j = 1, \dots, N_{\rm gr}$, in the graphene layer is reduced (circles in the inset box). Parameters: initial values $n^0 = p^0 = 100$, number of points $N = 820$, snapping parameter $\alpha = 2$, $\xi= 5\cdot 10^{-6}$, {\tt toll} = $10^{-5}$, $V_b = 0.2 V$ gate voltage $V_G = -1.5 V$ in panel (a) and $V_G = 1.5 V$ in panel (b).}}
   \label{fig:yprofile}
\end{figure}

To further investigate the effect of the discretization of the graphene layer, we consider the transverse profile of the electrostatic potential along the central section of the device,
\[
    \varphi\left(\frac{x_1+x_4}{2},y\right).
\]
Fig.~\ref{fig:yprofile} shows this profile for three different values of the effective thickness \(t_{\rm gr}\) of the graphene region. The total number of grid points \(N=820\), and thus also the space step $h$ are kept fixed, so that a reduction of \(t_{\rm gr}\) also reduces the number of grid points inside \(\Omega_{\rm gr}\). This provides a direct test of the sensitivity of the solution to the resolution of the active layer.

For negative gate voltage, the potential presents a local maximum in the graphene region, whereas for positive gate voltage, it presents a local minimum. The zoomed views highlight that the potential varies slightly across the graphene thickness and that this variation depends on \(t_{\rm gr}\). Therefore, even though the active layer is extremely thin, its effective thickness and vertical discretization influence the electrostatic profile. This supports the use of a full two-dimensional Poisson problem with two oxide/graphene interfaces, rather than a purely one-dimensional or lumped electrostatic approximation.

\section{Conclusions and Future work}
\label{SEC:concl}
In this work, we extended the ghost-FEM framework to elliptic interface problems with discontinuous diffusion coefficients and applied it to the numerical simulation of low-dimensional semiconductor devices. The proposed method combines an unfitted Cartesian-grid discretization with a level-set description of the geometry and a snapping-back-to-grid procedure to avoid ill-conditioning due to small cut cells. Boundary and interface conditions are imposed through a Nitsche-type variational formulation, allowing the method to handle internal interfaces without requiring a mesh fitted to the geometry.

The accuracy of the method was assessed through several benchmark elliptic interface problems involving different geometries and coefficient jumps. The numerical results show second-order convergence for the solution and first-order convergence for its gradient, confirming the expected accuracy of the spatial discretization. These tests demonstrate that the proposed approach provides an effective and flexible numerical framework for interface problems arising in heterogeneous media.

The method was then applied to a graphene field-effect transistor modeled by a self-consistent drift-diffusion--Poisson system. The electrostatic potential was computed in a two-dimensional oxide/graphene/oxide structure, while charge transport in the graphene layer was described by one-dimensional bipolar drift-diffusion equations. This setting naturally leads to a coupled problem involving different spatial dimensions, for which the ghost-FEM formulation and the domain-decomposition strategy provide a convenient computational framework.

The GFET simulations show that the gate voltage strongly affects the electrostatic potential and the resulting current response. The computed transfer characteristics exhibit a clear transition from a low-conducting regime to an ON state, while the ON current increases with the drain-source bias. Moreover, the analysis of the electrostatic potential inside the graphene layer shows that, although the layer is very thin, the transverse variation of the potential is not completely negligible. The effective graphene thickness and its vertical discretization influence the numerical solution, confirming the relevance of a full two-dimensional Poisson description with explicit oxide/graphene interface conditions.

Future work will address several extensions. A time-dependent formulation of the drift-diffusion model would provide additional flexibility in the choice of temporal discretization. In particular, high-order IMEX schemes \cite{boscarino2024implicit} could be employed to treat the nonlinear transport terms efficiently while maintaining the second-order accuracy achieved by the spatial derivatives shown in Fig.~\ref{fig:accuracy_3layers}. Further developments will also focus on more efficient implementations and mesh-refinement strategies, with the aim of resolving thinner graphene regions and more pronounced multiscale effects. Finally, the unfitted discretization framework considered in this work can be applied to refined meshes without the need for remeshing. This provides a basis for the investigation of semiconductor devices with more complex geometries and multiscale effects, which will be addressed in future work. The development of more efficient computational implementations could further enable high-resolution simulations, allowing the investigation of thinner graphene regions and a more accurate description of the associated multiscale phenomena.

\section*{Appendix}
\begin{table}
	\centering  
	\begin{tabular}{|c|c|c|c|c|c|}
		\hline
		Symbol       & value                            & Symbol        & value  & Symbol    & value\\ 
		\hline\hline 
        $\mu_0$ & $0.4650$ $\mu$m$^2$/V ps & $n_\mathrm{ref}$ & $1.1 \cdot 10^5$ $\mu$m$^{-2}$  & $\sigma$ &  2.2 \\
        \hline
		 $\epsilon_0$ & $8.85\cdot 10^{-18}$ & $\epsilon_{\rm ox}$  &   $3.6\epsilon_0$    &    $\epsilon_{\rm gr}$  &  $3.3\epsilon_0$ \\
		\hline 
		$L$ & 100\, nm & $y_2-y_1 = y_4-y_3$ & 10\,nm & $V_{G_u} = V_{G_d}$ & $V_G$\\
		\hline
        $n_{\rm imp}$ & $2.5\cdot 10^3 \mu$m$^{-2}$ & $e$ & $0.1602\cdot 10^{-18}$ C & $\Delta\varepsilon$  & 0.25 eV \\
		\hline 
	\end{tabular}
	\caption{\textit{Parameters involved.}}
	\label{table_parameters}
\end{table}
\subsection*{Steps of the two interfaces problems}
\label{appendix_steps}
The domain is divided in three subdomains, i.e. $\Omega = \Omega_1\cup\Omega_{\rm gr} \cup \Omega_2$ (see Fig.~\ref{fig:domain}). Usually, in domain decomposition procedures an iterative approach is considered to solve  system~\eqref{eq:variational2D_Wh}.
Since one grid line of \(\Omega_{\rm gr}\) is associated with the interface condition at \(\Gamma_1\) and another one with the interface condition at \(\Gamma_2\), at least three additional grid lines are required to perform the interior computation within \(\Omega_{\rm gr}\). For this reason, the number of cells in \(\Omega_{\rm gr}\) is chosen so that the subdomain contains at least five grid lines in total.

Below, we describe the computational steps to compute. We start from the central layer, with the same interface conditions at $\Gamma_1$ and $\Gamma_2$, because we expect symmetry respect to the $y$-direction. 
\begin{itemize}
    \item[\textbf{step 1:}] We solve the system for $\varphi_{\rm gr} \in \Omega_{\rm gr}$, with Dirichlet boundary conditions 
\begin{align}
\begin{split}        
-\nabla \cdot (\epsilon_{\rm gr} \nabla \varphi_{\rm gr}) &=e(p - n + n_{\rm imp})/t_{\rm gr} \qquad {\rm in }\,\Omega_{\rm gr} \\
        \varphi_{\rm gr} &= 0 \quad {\rm at} \, x = x_1 \\
        \varphi_{\rm gr} &= V_{\rm b} \quad {\rm at} \, x = x_4 \\ 
     \varphi_{\rm gr} & = 0  \quad {\rm at} \, \Gamma_1, \, \Gamma_2
\end{split}
\end{align}
    \item[\textbf{step 2:}]  We solve the system for $\varphi_1 \in \Omega_1$ and for $\varphi_2 \in \Omega_2$
\begin{align}
    \begin{split}
         -\nabla \cdot (\epsilon_{\rm ox} \nabla \varphi_i) &=0 \qquad {\rm in }\,\Omega_i \\
        \varphi_i &= 0 \quad {\rm at} \, x = x_1 \\
        \varphi_i &= V_{\rm b} \quad {\rm at} \, x = x_4 \\ 
        \varphi_1 & = V_{\rm G_d}  \quad {\rm at} \, y = y_1, \, x\in [x_2,x_3]  \\ 
         \varphi_2 & = V_{\rm G_u}  \quad {\rm at} \, y = y_4, \, x\in [x_2,x_3]  \\ 
       \epsilon_{\rm ox} \nabla_{\bf n} \varphi_i &  = \epsilon_{\rm gr} \nabla_{\bf n} \varphi_{\rm gr}   \quad {\rm at} \, \Gamma_i,   \\
        \nabla_{\bf n} \varphi_i & = 0 \quad {\rm in} \, \Gamma_{\rm N},
    \end{split}
\end{align}
    for $i = 1,2$.
    \item[\textbf{step 3:}]  We solve the system for $\varphi_{\rm gr} \in \Omega_{\rm gr}$, with Dirichlet boundary conditions 
\begin{align}
\begin{split}        
-\nabla \cdot (\epsilon_{\rm gr} \nabla \varphi_{\rm gr}) &=e(p - n + n_{\rm imp})/t_{\rm gr} \qquad {\rm in }\,\Omega_{\rm gr} \\
        \varphi_{\rm gr} &= 0 \quad {\rm at} \, x = x_1 \\
        \varphi_{\rm gr} &= V_{\rm b} \quad {\rm at} \, x = x_4 \\ 
     \varphi_{\rm gr} & = \varphi_i  \quad {\rm at} \, \Gamma_i, \quad i = 1,2
\end{split}
\end{align}
and repeat from \textbf{step 2}.
\end{itemize}

\section*{Acknowledgements}
The work of C.A. has been supported by the Spoke 10 Future AI Research (FAIR) of the Italian Research Center funded by the Ministry of University and Research as part of the National Recovery and Resilience Plan (PNRR). \\ 
C.A. has also been supported also by Italian Ministerial grant PRIN 2022 ``Efficient numerical schemes and optimal control methods for time-dependent partial differential equations''. %, No. 2022N9BM3N - Finanziato dall'Unione europea - Next Generation EU - CUP: E53D23005830006. \\
C.A. has also been supported by the Italian Ministerial grant PRIN 2022 PNRR ``FIN4GEO: Forward and Inverse Numerical Modeling of hydrothermal systems in volcanic regions with application to geothermal energy exploitation''. %, No. P2022BNB97 - Finanziato dall'Unione europea - Next Generation EU – CUP: E53D23017960001. \\
C.A. is member of the Gruppo Nazionale Calcolo Scientifico-Istituto Nazionale di Alta Matematica (GNCS-INdAM).\\
G.N. acknowledges the support from the Gruppo Nazionale per la Fisica Matematica - Istituto Nazionale di Alta Matematica (GNFM-INdAM).

\bibliographystyle{plain}
\bibliography{bibliography.bib}

\end{document}